\documentclass{article}
\usepackage{amscd}      % to be able to use "CD" environment
\usepackage{amssymb}
\usepackage{amsmath}
\usepackage{xypic}      % to be able to use "diagram" environment of
\LaTeXdiagrams          % xypic to typeset commutative diagrams
\usepackage[all,v2]{xy}
\xyoption{2cell} \UseAllTwocells \xyoption{frame} \CompileMatrices
\allowdisplaybreaks[3]
\usepackage{theorem}

\newtheorem{prop}{Proposition}[section]
\newtheorem{lem}[prop]{Lemma}

\newtheorem{cor}[prop]{Corollary}
\newtheorem{thm}[prop]{Theorem}

\theorembodyfont{\upshape}

\newtheorem{defn}[prop]{Definition}

\newtheorem{example}{Example}

\newtheorem{rmk}{Remark}

\newenvironment{pf}{\begin{trivlist}\item[]{\sc Proof.}}%
            {\nolinebreak $\Box$ \end{trivlist}}

\newcommand{\noprint}[1]{}

\newcommand{\toto}{\rightrightarrows}

\newcommand{\ldiag}[1]%
       {\makebox[0cm]{${\scriptstyle#1}\downarrow\phantom{\scriptstyle#1}$}}
\newcommand{\ldiagup}[1]%
       {\makebox[0cm]{${\scriptstyle#1}\uparrow\phantom{\scriptstyle#1}$}}
\newcommand{\rdiag}[1]%
       {\makebox[0cm]{$\phantom{\scriptstyle#1}\downarrow{\scriptstyle#1}$}}
\newcommand{\sediagr}[1]%
       {\makebox[0cm]{$\phantom{\scriptstyle#1}\searrow{\scriptstyle#1}$}}
\newcommand{\nediagr}[1]%
       {\makebox[0cm]{$\phantom{\scriptstyle#1}\nearrow{\scriptstyle#1}$}}
\newcommand{\rdiagup}[1]%
       {\makebox[0cm]{$\phantom{\scriptstyle#1}\uparrow{\scriptstyle#1}$}}
\newcommand{\swdiag}[1]%
       {\makebox[0cm]{$\phantom{\scriptstyle#1}\swarrow{\scriptstyle#1}$}}
\newcommand{\sediag}[1]%
       {\makebox[0cm]{${\scriptstyle#1}\searrow\phantom{\scriptstyle#1}$}}
\newcommand{\nediag}[1]%
       {\makebox[0cm]{${\scriptstyle#1}\nearrow\phantom{\scriptstyle#1}$}}

\newcommand{\doublearrowstack}[2]%
                      {{{{\scriptstyle#1}\atop{\textstyle\longrightarrow}}\atop{{\textstyle\longrightarrow}\atop{\scriptstyle#2}}}}
\newcommand{\rightleftarrowstack}[2]%
                      {{{{\scriptstyle#1}\atop{\textstyle\longrightarrow}}\atop{{\textstyle\longleftarrow}\atop{\scriptstyle#2}}}}
\newcommand{\leftrightarrowstack}[2]%
                      {{{{\scriptstyle#1}\atop{\textstyle\longleftarrow}}\atop{{\textstyle\longrightarrow}\atop{\scriptstyle#2}}}}

\newcommand{\overtoparrow}%
{\makebox[0cm]{\beginpicture \setcoordinatesystem units
<.8cm,.4cm> point at 0 0 \setplotarea x from -3 to 3, y from 0 to
1 \setquadratic \plot -3 0 0 1 3 0 / \put{\vector(3,-1){0}}[Bl] at
3 0
\endpicture}}

\newcommand{\underbottomarrow}%
{\makebox[0cm]{\beginpicture \setcoordinatesystem units
<.8cm,.4cm> point at 0 0 \setplotarea x from -3 to 3, y from 0 to
1 \setquadratic \plot -3 1 0 0 3 1 / \put{\vector(3,1){0}}[Bl] at
3 1
\endpicture}}

\newcommand{\ses}[5]%
{0\longrightarrow#1\stackrel{#2}{ \longrightarrow}#3\stackrel{#4}{
\longrightarrow}#5\longrightarrow0}

\newcommand{\dt}[6]%
{#1\stackrel{#2}{longrightarrow}#3
\stackrel{#4}{\longrightarrow}#5 \stackrel{#6}{\longrightarrow}
#1[1]}

\newcommand{\cat}[1]%
{(\mbox{\rm #1})}

\def\Label#1{\label{#1}{\tt [#1]}\phantom{h}}
\def\Label{\label}

\title{The Orbifold Cohomology Ring of Simplicial Toric Stack Bundles}
\author{Yunfeng Jiang}
\begin{document}
\sloppy \maketitle
\begin{abstract}
We introduce extended toric Deligne-Mumford stacks. We use an
extended toric Deligne-Mumford stack to get the toric stack bundle
and compute its orbifold Chow ring. Finally we generalize one
result of Borisov, Chen and Smith so that the orbifold Chow ring
of the toric stack bundle and the Chow ring of its crepant
resolution are fibres of a flat family.
\end{abstract}

%%% ----------------------------------------------------------------------
\maketitle
%%% ----------------------------------------------------------------------

\section{Introduction}
The Chen-Ruan orbifold cohomology was constructed by the genus
zero and degree zero orbifold Gromov-Witten invariants of
Deligne-Mumford stacks, see [5],[6],[1].  In this paper we discuss
the case of a toric stack bundle. The orbifold Chow ring of the
general toric Deligne-Mumford stack was obtained by Borisov, Chen
and Smith [4]. See also [11] for the case of weighted projective
space. Let $\mathbf{\Sigma}:=(N,\Sigma,\beta)$ be a stacky fan.
The  toric Deligne-Mumford stack
$\mathcal{X}(\mathbf{\Sigma})=[Z/G]$ is a quotient stack, where
$Z=\mathbb{C}^{n}-V(J_{\Sigma})$ and $J_{\Sigma}$ is the
square-free ideal of $\Sigma$. The action of $G$  on $Z$ is
through the map $\alpha: G \longrightarrow
(\mathbb{C}^{\times})^{n}$ determined by the stacky fan.  Let
$P\longrightarrow B$ be a principal
$(\mathbb{C}^{\times})^{n}$-bundle over a smooth variety $B$,
define $^{P}\mathcal{X}(\mathbf{\Sigma})$ to be the quotient stack
$[(P\times_{(\mathbb{C}^{\times})^{n}}Z)/G]$, where $G$ acts on
$P$ trivially. The stack
$[(P\times_{(\mathbb{C}^{\times})^{n}}Z)/G]$ may be written as
$P\times_{(\mathbb{C}^{\times})^{n}}[Z/G]$, then
$^{P}\mathcal{X}(\mathbf{\Sigma})\longrightarrow B$ is a toric
stack bundle over $B$ with fibre the  toric Deligne-Mumford stack
$\mathcal{X}(\mathbf{\Sigma})$. We study the orbifold Chow ring of
$^{P}\mathcal{X}(\mathbf{\Sigma})$.

Before we go further, let's consider the case when $N$ is a finite
abelian group, then in the stacky fan
$\mathbf{\Sigma}=(N,\Sigma,\beta)$, $\Sigma=0$ and $\beta$ is the
zero homomorphism $0\longrightarrow N$. The toric Deligne-Mumford
stack $\mathcal{X}(\mathbf{\Sigma})=[pt/G]$ is the classifying
stack $\mathcal{B}G$, where
$G=Hom_{\mathbb{Z}}(N,\mathbb{C}^{\times})$.  As a stack
$\mathcal{B}G$ can have different representations. For example, if
$N=\mathbb{Z}/3\mathbb{Z}$, $G=\mu_{3}$ is the cyclic group of
order $3$, then $\mathcal{B}\mu_{3}$ can also be represented by
the quotient stack $[\mathbb{C}^{\times}/\mathbb{C}^{\times}]$,
where the action is given by $(\cdot)^{3}$. The stack
$[\mathbb{C}^{\times}/\mathbb{C}^{\times}]$ is not a toric
Deligne-Mumford stack in the sense of [4].  If
$\mathcal{B}\mu_{3}=[\mathbb{C}^{\times}/\mathbb{C}^{\times}]$,
let $L\longrightarrow B$ be a line bundle over a smooth variety
$B$, then
$[(L^{\times}\times_{\mathbb{C}^{\times}}\mathbb{C}^{\times})/\mathbb{C}^{\times}]$
is a nontrivial $\mu_{3}$-gerbe over $B$ if the line bundle is
nontrivial. While if $\mathcal{B}\mu_{3}=[pt/\mu_{3}]$, we can't
twist it by any line bundle.

In order to make
$\mathcal{B}\mu_{3}=[\mathbb{C}^{\times}/\mathbb{C}^{\times}]$ a
toric Deligne-Mumford stack,  we slightly generalize the
construction of toric Deligne-Mumford stacks. We introduce
$extended ~stacky~ fan~
\mathbf{\Sigma^{e}}:=(N,\Sigma,\beta^{e})$, where $N$ and $\Sigma$
are  the same as in a stacky fan
$\mathbf{\Sigma}=(N,\Sigma,\beta)$, but $\beta^{e}:
\mathbb{Z}^{m}\longrightarrow N$ is determined by
$b_{1},\ldots,b_{n},b_{n+1},\ldots,b_{m}\in N$ satisfying the
conditions that $m\geq n$, $\overline{b}_{i}$ generates the ray
$\rho_{i}$ for $1\leq i\leq n$ and all the other data
$\{b_{n+1},\ldots,b_{m}\}$ belong to $N$, where
$\overline{b}_{i}\in\overline{N}$, and  $\overline{N}$ is the
lattice in $N_{\mathbb{Q}}$ determined by the projection
$N\longrightarrow \overline{N}$.  We call
$\{b_{n+1},\cdots,b_{m}\}$ the extra data in
$\mathbf{\Sigma^{e}}$. Associated to an extended stacky fan
$\mathbf{\Sigma^{e}}$, we define the extended toric
Deligne-Mumford stack
$\mathcal{X}(\mathbf{\Sigma^{e}}):=[Z^{e}/G^{e}]$ as a quotient
stack, where $Z^{e}=Z\times (\mathbb{C}^{\times})^{m-n}$ and
$G^{e}$ acts on $Z^{e}$ through the homomorphism $\alpha^{e}:
G^{e}\longrightarrow (\mathbb{C}^{\times})^{m}$  determined by the
extended stacky fan. It is easy to see that every extended stacky
fan $\mathbf{\Sigma^{e}}$ naturally determines a stacky fan
$\mathbf{\Sigma}$. We prove that the extended toric
Deligne-Mumford stack $\mathcal{X}(\mathbf{\Sigma^{e}})$ is
isomorphic to the underlying toric Deligne-Mumford stack
$\mathcal{X}(\mathbf{\Sigma})$. But we have more freedom to twist
as our example above shows.

Given an extended toric Deligne-Mumford stack
$\mathcal{X}(\mathbf{\Sigma^{e}})$, from the extended stacky fan
$\mathbf{\Sigma^{e}}$, we have the following exact sequence:
$$1\longrightarrow \mu\longrightarrow
G^{e}\stackrel{\alpha^{e}}{\longrightarrow}
(\mathbb{C}^{\times})^{m}\longrightarrow T\longrightarrow 1$$
where $T=(\mathbb{C}^{\times})^{d}$. Let $P\longrightarrow B$ be a
principal $(\mathbb{C}^{\times})^{m}$-bundle,  let
$^{P}\mathcal{X}(\mathbf{\Sigma^{e}})$  be the quotient stack
$[(P\times_{(\mathbb{C}^{\times})^{m}}Z^{e})/G^{e}]$, where
$G^{e}$ acts on $P$ trivially and on $(\mathbb{C}^{\times})^{m}$
through the map $\alpha^{e}$ in the above exact sequence.  Then
$^{P}\mathcal{X}(\mathbf{\Sigma^{e}})$ is a toric stack bundle
over $B$ with fibre the extended toric Deligne-Mumford stack
$\mathcal{X}(\mathbf{\Sigma^{e}})$.  The extra data
$\{b_{n+1},\cdots,b_{m}\}$ in $\mathbf{\Sigma^{e}}$ can be put
into the $Box(\mathbf{\Sigma^{e}})$ which do not influence the
structure of the toric stack bundle
$^{P}\mathcal{X}(\mathbf{\Sigma^{e}})$. The choice of torsion and
nontorsion extra data does affect the structure of
$^{P}\mathcal{X}(\mathbf{\Sigma^{e}})$, but does not affect the
orbifold cohomology.  To describe the orbifold Chow ring of
$^{P}\mathcal{X}(\mathbf{\Sigma^{e}})$, we introduce some line
bundles over $B$. Let $M=N^{*}$ be the dual of $N$. For $\theta\in
M$, let $\xi_{\theta}\longrightarrow B$ be the line bundle coming
from the principal $T$ bundle $E\longrightarrow B$ by "extending"
the structure group via $\chi^{\theta}:T\longrightarrow
\mathbb{C}^{\times}$, where $E\longrightarrow B$ is induced from
the $(\mathbb{C}^{\times})^{m}$-bundle $P$ in the above exact
sequence.  Define the deformed ring
$A^{*}(B)[N]^{\mathbf{\Sigma^{e}}}=A^{*}(B)\otimes
\mathbb{Q}[N]^{\mathbf{\Sigma^{e}}}$ where
$\mathbb{Q}[N]^{\mathbf{\Sigma^{e}}}:=\bigoplus_{c\in N}y^{c}$,
$y$ is the formal variable and $A^{*}(B)$ is the Chow ring of $B$.
The multiplication of $\mathbb{Q}[N]^{\mathbf{\Sigma^{e}}}$ is
given by:
\begin{equation}\Label{product}
y^{c_{1}}\cdot y^{c_{2}}:=\begin{cases}y^{c_{1}+c_{2}}&\text{if
there is a cone}~ \sigma\in\Sigma ~\text{such that}~ \overline{c}_{1}\in\sigma, \overline{c}_{2}\in\sigma\,,\\
0&\text{otherwise}\,.\end{cases}
\end{equation}
Let $\mathcal{I}(\mathbf{\Sigma^{e}})$ be the ideal in
$A^{*}(B)[N]^{\mathbf{\Sigma^{e}}}$ generated by the elements:
\begin{equation}\Label{ideal}
\left(c_{1}(\xi_{\theta})+\sum_{i=1}^{n}\theta(b_{i})y^{b_{i}}\right)_{\theta\in
M}
\end{equation} for $\theta\in M$ and
$A_{orb}^{*}\left(^{P}\mathcal{X}(\mathbf{\Sigma^{e}})\right)$ be
the orbifold Chow ring of the toric stack bundle. Then we have the
following Theorem:

\begin{thm}
If $^{P}\mathcal{X}(\mathbf{\Sigma^{e}})\longrightarrow B$ is a
toric stack bundle over a smooth variety $B$ with fibre the
extended  toric Deligne-Mumford stack
$\mathcal{X}(\mathbf{\Sigma^{e}})$ associated to an extended
stacky fan $\mathbf{\Sigma^{e}}$, then we have an isomorphism of
$\mathbb{Q}$-graded rings:
$$A_{orb}^{*}\left(^{P}\mathcal{X}(\mathbf{\Sigma^{e}})\right)\cong \frac{A^{*}(B)[N]^{\mathbf{\Sigma^{e}}}}{\mathcal{I}(\mathbf{\Sigma^{e}})}$$
\end{thm}

To prove this theorem, first using the similar result in [4] that
the components of the inertia stack
$\mathcal{I}(\mathcal{X}(\mathbf{\Sigma^{e}}))$ of
$\mathcal{X}(\mathbf{\Sigma^{e}})$ is given by
$Box(\mathbf{\Sigma^{e}})$ which determines all the elements in
the local group of $\mathcal{X}(\mathbf{\Sigma^{e}})$, we explain
that the twist by the $(\mathbb{C}^{\times})^{m}$-bundle $P$ does
not twist the components of the inertia stack of the toric stack
bundle $^{P}\mathcal{X}(\mathbf{\Sigma^{e}})$. Then this makes it
possible to use the similar methods as in [4] to determine
3-twisted sectors, obstruction bundles of
$^{P}\mathcal{X}(\mathbf{\Sigma^{e}})$ and compute the orbifold
Chow ring of $^{P}\mathcal{X}(\mathbf{\Sigma^{e}})$. As an
example, let $N$ be a finite abelian group and
$\mathbb{Z}\longrightarrow N$ be any homomorphism. Then
$\mathbf{\Sigma^{e}}=(N,0,\beta^{e})$ is an extended stacky fan,
and $\mathcal{X}(\mathbf{\Sigma^{e}})=\mathcal{B}\mu$, where
$\mu=Hom(N,\mathbb{C}^{\times})$. Twist  this extended toric
Deligne-Mumford stack by a line bundle $L$ over a smooth variety
$B$, we get the $\mu$-gerbe over $B$. We determine its inertia
stack and compute its orbifold Chow ring.

The paper is organized  as follows. In Section 2 we introduce
extended toric Deligne-Mumford stacks. In Section 3 we define the
toric stack bundle and discuss its properties. In Section 4 we
describe the orbifold Chow ring of the toric stack bundle. In
Section 5 we give an interesting example of toric stack bundle,
the $\mu$-gerbe $\mathcal{X}$ over $B$ for a finite abelian group
$\mu$ and smooth variety $B$. Finally in Section 6 we give some
applications of crepant resolutions.

In this paper, we use the rational numbers $\mathbb{Q}$ as
coefficients of the Chow ring and orbifold Chow ring.  By an
$orbifold$ we mean a smooth Deligne-Mumford stack such that at the
generic point the automorphism group is trivial. This type of
orbifold is sometimes called  a $reduced~ orbifold$ in
differential geometry.

%%% ----------------------------------------------------------------------
\subsection*{Acknowledgments}
I would like to thank my  advisor Kai Behrend for his help in
preparing  this work. I also thank Andrei Mustata and Hsian-Hua
Tseng for valuable discussions.
%%% ----------------------------------------------------------------------
\section{The Extended Toric Deligne-Mumford stacks.}

In this section we introduce extended stacky fans and construct
extended toric Deligne-Mumford stacks. We prove that the extended
toric Deligne-Mumford stack is isomorphic to the underlying  toric
Deligne-Mumford stack.

We refer to [4] the construction and notation of toric
Deligne-Mumford stacks.  Let $N$ be a finitely generated abelian
group of $rank~ d$. Let $\overline{N}$ be the lattice generated by
$N$ in the $d$-dimensional vector space
$N_{\mathbb{Q}}:=N\otimes_{\mathbb{Z}}\mathbb{Q}$. The natural map
$N\longrightarrow \overline{N}$ is denoted by $b\longmapsto
\overline{b}$. Let $\Sigma$ be a rational $simplicial$ fan in
$N_{\mathbb{Q}}$. Suppose  $\rho_{1},\ldots,\rho_{n}$ are the rays
in $\Sigma$. We fix $b_{i}\in N$ for $1\leq i\leq n$ such that
$\overline{b}_{i}$ generates the cone $\rho_{i}$. We  choose extra
data $\{b_{n+1},\ldots,b_{m}\}\subset N$ and  consider the
homomorphism $\beta^{e}: \mathbb{Z}^{m}\longrightarrow N$
determined by the elements $b_{1},\ldots,b_{m}$. We require  that
$\beta^{e}$ has finite cokernel.
$\mathbf{\Sigma^{e}}:=(N,\Sigma,\beta^{e})$ is called an $extended
~stacky~ fan$.

It is easy to see that any extended stacky fan
$\mathbf{\Sigma^{e}}=(N,\Sigma,\beta^{e})$ naturally determines a
stacky fan $\mathbf{\Sigma}:=(N,\Sigma,\beta)$, where $\beta:
\mathbb{Z}^{n}\longrightarrow N$ is given by
$b_{1},\ldots,b_{n}\in N$. Now since $\beta^{e}$ has finite
cokernel, from Proposition 2.2 in [4], we have exact sequences:
$$0\longrightarrow DG(\beta^{e})^{*}\longrightarrow \mathbb{Z}^{m}\stackrel{\beta^{e}}{\longrightarrow} N\longrightarrow Coker(\beta^{e})\longrightarrow
0$$
$$0\longrightarrow N^{*}\longrightarrow \mathbb{Z}^{m}\stackrel{(\beta^{e})^{\vee}}{\longrightarrow} DG(\beta^{e})\longrightarrow Coker((\beta^{e})^{\vee})\longrightarrow 0$$
where $(\beta^{e})^{\vee}$ is the Gale dual of $\beta^{e}$. As a
$\mathbb{Z}$-module, $\mathbb{C}^{\times}$ is divisible, so it is
an injective $\mathbb{Z}$-module, from [13], the functor
$Hom_{\mathbb{Z}}(-,\mathbb{C}^{\times})$ is exact.  We get the
exact sequence:
$$1\longrightarrow Hom_{\mathbb{Z}}(Coker((\beta^{e})^{\vee}),\mathbb{C}^{\times})\longrightarrow
Hom_{\mathbb{Z}}(DG(\beta^{e}),\mathbb{C}^{\times})\longrightarrow
Hom_{\mathbb{Z}}(\mathbb{Z}^{m},\mathbb{C}^{\times})$$
$$\longrightarrow
Hom_{\mathbb{Z}}(N^{*},\mathbb{C}^{\times})\longrightarrow 1$$ Let
$\mu:=Hom_{\mathbb{Z}}(Coker((\beta^{e})^{\vee}),\mathbb{C}^{\times})$,
we have the exact sequence:
\begin{equation}\Label{stack}
1\longrightarrow \mu\longrightarrow
G^{e}\stackrel{\alpha^{e}}{\longrightarrow}
(\mathbb{C}^{\times})^{m}\longrightarrow T\longrightarrow 1 \
\end{equation}
From [4], the toric Deligne-Mumford stack
$\mathcal{X}(\mathbf{\Sigma})=[Z/G]$ is a quotient stack, where
they use the method of quotient construction of toric varieties
[7]. Define $Z^{e}:=Z\times (\mathbb{C}^{\times})^{m-n}$, then
there exists a natural action of $(\mathbb{C}^{\times})^{m}$ on
$Z^{e}$. The group $G^{e}$ acts on $Z^{e}$ through the map
$\alpha^{e}$ in (\ref{stack}). The quotient stack $[Z^{e}/G^{e}]$
is associated to the groupoid $Z^{e}\times G^{e}\toto Z^{e}$.
Define the morphism $\varphi: Z^{e}\times G^{e}\longrightarrow
Z^{e}\times Z^{e}$ to be $\varphi(x,g)=(x,g\cdot x)$. Since
$Z^{e}=Z\times (\mathbb{C}^{\times})^{m-n}$, we can mimic the
proof the   Lemma 3.1 in [4] to get that $\varphi$ is finite. So
the stack $[Z^{e}/G^{e}]$ is a Deligne-Mumford stack.
\begin{lem}
The morphism $\varphi: Z^{e}\times G^{e}\longrightarrow
Z^{e}\times Z^{e}$ is a finite morphism.
\end{lem}

\begin{defn}
For an extended stacky fan
$\mathbf{\Sigma^{e}}=(N,\Sigma,\beta^{e})$, define the extended
toric Deligne-Mumford stack $\mathcal{X}(\mathbf{\Sigma^{e}})$ to
be the quotient stack $[Z^{e}/G^{e}]$.
\end{defn}

\begin{prop}
For an extended stacky fan
$\mathbf{\Sigma^{e}}=(N,\Sigma,\beta^{e})$, the extended toric
Deligne-Mumford stack $\mathcal{X}(\mathbf{\Sigma^{e}})$ is
isomorphic to the underlying toric Deligne-Mumford stack
$\mathcal{X}(\mathbf{\Sigma})$.
\end{prop}
\begin{pf}
From the definitions of extended stacky fan $\mathbf{\Sigma^{e}}$
and stacky fan $\mathbf{\Sigma}$, we have the following
commutative diagram:
\[
\begin{CD}
0 @ >>>\mathbb{Z}^{n}@ >>> \mathbb{Z}^{m}@ >>> \mathbb{Z}^{m-n} @
>>> 0\\
&& @VV{\beta}V@VV{\beta^{e}}V@VV{\widetilde{\beta}}V \\
0@ >>> N @ >{id}>>N@ >>> 0 @>>> 0
\end{CD}
\]
From the definition of Gale dual, we compute
$DG(\widetilde{\beta})=\mathbb{Z}^{m-n}$ and
$\widetilde{\beta}^{\vee}$ is an isomorphism. So from   Lemma 2.3
in [4], applying the Gale dual and the
$Hom_{\mathbb{Z}}(-,\mathbb{C}^{\times})$ functor to the above
diagram we get:
\begin{equation}\Label{second}
\begin{CD}
1 @ >>>G@ >{\varphi_{1}}>> G^{e}@ >>> (\mathbb{C}^{\times})^{m-n}
@
>>> 1\\
&& @VV{\alpha}V@VV{\alpha^{e}}V@VV{\widetilde{\alpha}}V \\
1@ >>>(\mathbb{C}^{\times})^{n} @ >>>(\mathbb{C}^{\times})^{m}@
>>>(\mathbb{C}^{\times})^{m-n} @>>> 1
\end{CD}
\end{equation}
We define the morphism $\varphi_{0}: Z\longrightarrow
Z^{e}=Z\times (\mathbb{C}^{\times})^{m-n}$ to be the inclusion
defined by $z\longmapsto (z,1)$. So $(\varphi_{0}\times
\varphi_{1}, \varphi_{0}): (Z\times G\toto Z)\longrightarrow
(Z^{e}\times G^{e}\toto Z^{e})$ defines a morphism between
groupoids.  Let $\varphi: [Z/G]\longrightarrow [Z^{e}/G^{e}]$ be
the morphism of stacks induced from $(\varphi_{0}\times
\varphi_{1}, \varphi_{0})$. From the above commutative diagram we
have the following commutative diagram:
$$\xymatrix{
Z\times G\dto^{(s,t)}\rto^{\varphi_{0}\times\varphi_{1}} &Z^{e}\times G^{e}\dto^{(s,t)}\\
Z\times Z\rto^{\varphi_{0}\times\varphi_{0}} &Z^{e}\times Z^{e}
}$$ In (\ref{second}), $\widetilde{\alpha}$ is an isomorphism
which implies that the left square in (\ref{second}) is cartesian.
So the above commutative diagram is cartesian. $\varphi:
[Z/G]\longrightarrow [Z^{e}/G^{e}]$ is injective. Given an element
$(z_{1},\ldots,z_{n},z_{n+1},\ldots,z_{m})\in Z^{e}$, there exists
an element $g^{e}\in (\mathbb{C}^{\times})^{m-n}$ such that
$g^{e}\cdot(z_{1},\ldots,z_{n},z_{n+1},\ldots,z_{m})=(z_{1},\ldots,z_{n},1,\ldots,1)$.
From (\ref{second}), $g^{e}$ determine an element in $G^{e}$, so
$\varphi$ is surjective.  The  stacks
$\mathcal{X}(\mathbf{\Sigma^{e}})$ and
$\mathcal{X}(\mathbf{\Sigma})$ are isomorphic.
\end{pf}

Let $X(\Sigma)$ be the simplicial toric variety associated to the
extended stacky fan $\mathbf{\Sigma^{e}}$. We have the following
corollaries:
\begin{cor}
Given an extended stacky fan $\mathbf{\Sigma^{e}}$, then the
coarse moduli space of the extended toric Deligne-Mumford stack
$\mathcal{X}(\mathbf{\Sigma^{e}})$ is also the simplicial toric
variety $X(\Sigma)$.
\end{cor}

\begin{rmk}
As in [4], let $\sigma$ be a top dimensional cone in $\Sigma$,
denote by $Box(\sigma)$ to be the set of elements $v\in N$ such
that $\overline{v}=\sum_{\rho_{i}\subseteq
\sigma}a_{i}\overline{b}_{i}$ for some $0\leq a_{i}<1$. The set
$Box(\sigma)$ is in one-to-one correspondence with the elements in
the finite group $N(\sigma)=N/N_{\sigma}$, where $N(\sigma)$ is a
local group of the stack $\mathcal{X}(\mathbf{\Sigma^{e}})$. If
$\tau\subseteq \sigma$ is a low dimensional cone, we define
$Box(\tau)$ to be the set of elements in $v\in N$ such that
$\overline{v}=\sum_{\rho_{i}\subseteq \tau}a_{i}\overline{b}_{i}$,
where $0\leq a_{i}<1$. It is easy to see that $Box(\tau)\subset
Box(\sigma)$. In fact the elements in $Box(\tau)$ generate a
subgroup of the local group $N(\sigma)$. Let
$Box(\mathbf{\Sigma^{e}})$ be the union of $Box(\sigma)$ for all
$d$-dimensional cones $\sigma\in \Sigma$. For
$v_{1},\ldots,v_{n}\in N$, let
$\sigma(\overline{v}_{1},\ldots,\overline{v}_{n})$ be the unique
minimal cone in $\Sigma$ containing
$\overline{v}_{1},\ldots,\overline{v}_{n}$.
\end{rmk}

%%% --------------------------------------------------------------------------
%%%----------------------------------------------------------------------------
\section{The Toric Stack Bundle $^{P}\mathcal{X}(\mathbf{\Sigma^{e}})$.}

In this section we introduce the toric stack bundle
$^{P}\mathcal{X}(\mathbf{\Sigma^{e}})$ and determine its twisted
sectors.  Let $P\longrightarrow B$ be a principal
$(\mathbb{C}^{\times})^{m}$-bundle over a smooth variety $B$.  We
give the following definition.
\begin{defn}
We define the toric stack bundle
$^{P}\mathcal{X}(\mathbf{\Sigma^{e}})\longrightarrow B$ to be the
quotient stack
\begin{equation}\Label{bundle}
^{P}\mathcal{X}(\mathbf{\Sigma^{e}}):=[(P\times_{(\mathbb{C}^{\times})^{m}}Z^{e})/G^{e}]
\end{equation}
where $G^{e}$ acts on $P$ trivially.
\end{defn}

\begin{rmk}
Let $\phi: \mathbb{Z}^{m}\longrightarrow \mathbb{Z}^{m}$ be the
map given by $e_{i}\longmapsto e_{i}$ for $1\leq i\leq n$, and
$e_{j}\longmapsto e_{j}+\sum_{i=1}^{n}a_{i}^{j}e_{i}$ for $n+1\leq
j\leq m$, where $a_{i}^{j}\in \mathbb{Z}$. Then consider the
following commutative diagram:
\[
\begin{CD}
0 @ >>>\mathbb{Z}^{m}@ >{\phi}>> \mathbb{Z}^{m}@ >>> 0 @
>>> 0\\
&& @VV{\widetilde{\beta}^{e}}V@VV{\beta^{e}}V@VV{id}V \\
0@ >>> N @ >{id}>>N@ >>> 0 @>>> 0
\end{CD}
\]
We obtain a new extended stacky fan
$\mathbf{\widetilde{\Sigma}^{e}}=(N,\Sigma,\widetilde{\beta}^{e})$,
where the extra data in $\mathbf{\widetilde{\Sigma}^{e}}$ are
$b_{n+1}^{'}=b_{n+1}+\sum_{i=1}^{n}a_{i}^{n+1}b_{i},\cdots,b_{m}^{'}=b_{m}+\sum_{i=1}^{n}a_{i}^{m}b_{i}$.
 The map $\phi$ gives a map
$\mathbb{C}^{n}\times (\mathbb{C}^{\times})^{m-n}\longrightarrow
\mathbb{C}^{n}\times (\mathbb{C}^{\times})^{m-n}$ which is the
identity on the first factor and given by $\phi$ on the second
factor. Since the map in the above diagram doesn't change the fan
in the extended stacky fans, we have a map $\varphi_{0}:
P\times_{(\mathbb{C}^{\times})^{m}}Z^{e}\longrightarrow
P\times_{(\mathbb{C}^{\times})^{m}}Z^{e}$, we use the same proof
in Proposition 2.3 to prove that
$^{P}\mathcal{X}(\mathbf{\Sigma^{e}})\cong
~^{P}\mathcal{X}(\mathbf{\widetilde{\Sigma}^{e}})$. So this means
that we always can choose the extra data such that
$b_{j}=\sum_{i=1}^{n}a_{i}b_{i}$ for $j=n+1,\cdots,m$ and $0\leq
a_{i}<1$. These extra data are actually  in the
$Box(\mathbf{\Sigma^{e}})$.
\end{rmk}

\begin{example}
From the above Remark, the extra data can be put into
$Box(\mathbf{\Sigma^{e}})$. In this example we prove that they can
not be put into the torsion subgroup of $N$.  Let $N=\mathbb{Z}$
and $b_{1}=2, b_{2}=-2$. Then $\Sigma=\{b_{1},b_{2}\}$ is a
simplicial fan in $N_{\mathbb{Q}}$. Let
$\mathbf{\Sigma^{e}}=(N,\Sigma,\beta^{e})$, where $\beta^{e}:
\mathbb{Z}^{3}\longrightarrow \mathbb{Z}$ is determined by
$\{b_{1},b_{2},b_{3}=1\}$, then we compute
$DG(\beta^{e})=\mathbb{Z}^{2}$ and the Gale dual
$(\beta^{e})^{\vee}: \mathbb{Z}^{3}\longrightarrow \mathbb{Z}^{2}$
is given by the matrix $\left[
\begin{array}{ccc}
1&1&0\\
-1&0&2
\end{array}
\right]$. From Section 2, the extended toric Deligne-Mumford stack
$\mathcal{X}(\mathbf{\Sigma^{e}})=[(\mathbb{C}^{2}-\{0\})\times
\mathbb{C}^{\times}/(\mathbb{C}^{\times})^{2}]$, where the action
is given by
$(\lambda_{1},\lambda_{2})(x,y,z)=(\lambda_{1}\lambda_{2}^{-1}\cdot
x,\lambda_{1}\cdot y,\lambda_{2}^{2}\cdot z)$. We get
$\mathcal{X}(\mathbf{\Sigma^{e}})=\mathbb{P}^{1}\times[\mathbb{C}^{\times}/\mathbb{C}^{\times}]=
\mathbb{P}^{1}\times \mathcal{B}\mu_{2}$.  Now let
$\mathbf{\widetilde{\Sigma}^{e}}=(N,\Sigma,\widetilde{\beta}^{e})$,
where $\widetilde{\beta}^{e}: \mathbb{Z}^{3}\longrightarrow
\mathbb{Z}$ is determined by
$\{b_{1},b_{2},\widetilde{b}_{3}=0\}$, then we compute
$DG(\beta^{e})=\mathbb{Z}^{2}\oplus \mathbb{Z}_{2}$ and the Gale
dual $(\widetilde{\beta}^{e})^{\vee}:
\mathbb{Z}^{3}\longrightarrow \mathbb{Z}^{2}\oplus \mathbb{Z}_{2}$
is given by the matrix $\left[
\begin{array}{ccc}
1&1&0\\
0&0&1\\
0&0&0
\end{array}
\right]$.  The extended toric Deligne-Mumford stack
$\mathcal{X}(\mathbf{\widetilde{\Sigma}^{e}})=[(\mathbb{C}^{2}-\{0\})\times
\mathbb{C}^{\times}/(\mathbb{C}^{\times})^{2}\times \mu_{2}]$,
where the action is given by
$(\lambda_{1},\lambda_{2},\lambda_{3})(x,y,z)=(\lambda_{1}\cdot
x,\lambda_{1}\cdot y,\lambda_{2}\cdot z)$. We get
$\mathcal{X}(\mathbf{\widetilde{\Sigma}^{e}})=[\mathbb{P}^{1}/\mu_{2}]=
\mathbb{P}^{1}\times \mathcal{B}\mu_{2}$. Let $B=\mathbb{P}^{1}$
and $P=\mathbb{C}^{\times}\oplus \mathbb{C}^{\times}\oplus
\mathcal{O}(-1)^{\times}$, then
$^{P}\mathcal{X}(\mathbf{\Sigma^{e}})$ is a nontrivial
$\mu_{2}$-gerbe over $\mathbb{P}^{1}\times \mathbb{P}^{1}$ coming
from the line bundle $\mathcal{O}_{\mathbb{P}^{1}\times
\mathbb{P}^{1}}(0,-1)$. Let $Q=\mathcal{O}(n_{1})^{\times}\oplus
\mathcal{O}(n_{2})^{\times}\oplus \mathcal{O}(n_{3})^{\times}$,
then  $^{Q}\mathcal{X}(\mathbf{\widetilde{\Sigma}^{e}})$ is the
trivial $\mu_{2}$-gerbe over the  $\mathbb{P}^{1}$-bundle $E$ over
$\mathbb{P}^{1}$. So $^{P}\mathcal{X}(\mathbf{\Sigma^{e}})$ is not
isomorphic to $^{Q}\mathcal{X}(\mathbf{\widetilde{\Sigma}^{e}})$
for any $Q$.
\end{example}

From Corollary 2.4, $\mathcal{X}(\mathbf{\Sigma^{e}})$ has the
coarse moduli space $X(\Sigma)$ which is the simplicial toric
variety associated to the simplicial fan $\Sigma$.  From the exact
sequence in (\ref{stack}),  a $(\mathbb{C}^{\times})^{m}$-bundle
over $B$ determine a $T$-bundle over $B$ naturally. Let
$E\longrightarrow B$ be the principal $T$-bundle induced by $P$,
then we have the twists
$^{P}\mathcal{X}_{red}(\mathbf{\Sigma^{e}})\longrightarrow B$ with
fibre the  toric orbifold $\mathcal{X}_{red}(\mathbf{\Sigma^{e}})$
and $^{E}X(\Sigma)\longrightarrow B$ with fibre the simplicial
toric variety $X(\Sigma)$, where
$^{P}\mathcal{X}_{red}(\mathbf{\Sigma^{e}}):=[(P\times_{(\mathbb{C}^{\times})^{m}}Z^{e})/\overline{G}^{e}]$
and $^{E}X(\Sigma):=E\times_{T}X(\Sigma)$, and
$\overline{G}^{e}=Im(\alpha^{e})$ in (\ref{stack}), so we obtain
the exact sequence:
\begin{equation}\Label{sequence}
1\longrightarrow \mu\longrightarrow
G^{e}\stackrel{\alpha^{e}}{\longrightarrow} \overline{G}^{e}
\longrightarrow 1 \
\end{equation}
From [8], we have:
\begin{prop} $^{P}\mathcal{X}(\mathbf{\Sigma^{e}})$ is a
$\mu$-gerbe over $^{P}\mathcal{X}_{red}(\mathbf{\Sigma^{e}})$ for
a finite abelian group $\mu$.
\end{prop}

\begin{rmk}
In fact, any extended toric Deligne-Mumford stack is $\mu$-gerbe
over the underlying toric orbifold for a finite abelian group
$\mu$ and some kind of $\mu$-gerbes over toric Deligne-Mumford
stacks are again toric Deligne-Mumford stacks, see [3].
\end{rmk}

Because any toric stack bundle is a $\mu$-gerbe over the
corresponding toric orbifold bundle and can be represented as a
quotient stack, we have the following propositions:
\begin{prop}
The simplicial toric bundle $^{E}X(\Sigma)$ is the coarse moduli
space of the toric stack bundle
$^{P}\mathcal{X}(\mathbf{\Sigma^{e}})$ and the toric orbifold
bundle $^{P}\mathcal{X}_{red}(\mathbf{\Sigma^{e}})$.
\end{prop}
\begin{pf}
The toric stack bundle $^{P}\mathcal{X}(\mathbf{\Sigma^{e}})$ is a
$\mu$-gerbe over the simplicial toric orbifold bundle
$^{P}\mathcal{X}_{red}(\mathbf{\Sigma^{e}})$ for a finite abelian
group $\mu$,
$^{P}\mathcal{X}(\mathbf{\Sigma^{e}})=[(P\times_{(\mathbb{C}^{\times})^{m}}Z^{e})/G^{e}]$
and
$^{P}\mathcal{X}_{red}(\mathbf{\Sigma^{e}})=[(P\times_{(\mathbb{C}^{\times})^{m}}Z^{e})/\overline{G}^{e}]$
are quotient stacks. Take the  geometric quotient, we have the
coarse moduli space
$(P\times_{(\mathbb{C}^{\times})^{m}}Z^{e})//\overline{G}^{e}=(P\times
Z^{e})//(\mathbb{C}^{\times})^{m}\times \overline{G}^{e}$. From
Proposition 2.1 in [3], we have
$X(\Sigma)=Z//\overline{G}=Z^{e}//\overline{G}^{e}$, so
$$
E\times_{T}(Z^{e}//\overline{G}^{e})
=(P\times_{(\mathbb{C}^{\times})^{m}}T)\times_{T}(Z^{e}//\overline{G}^{e})
=(P\times Z^{e})//(\mathbb{C}^{\times})^{m}\times \overline{G}^{e}
$$
From the universal geometric quotients in [14], $^{E}X(\Sigma)$ is
the coarse moduli space of $^{P}\mathcal{X}(\mathbf{\Sigma^{e}})$
and $^{P}\mathcal{X}_{red}(\mathbf{\Sigma^{e}})$.
\end{pf}

\begin{prop}
The  toric stack bundle $^{P}\mathcal{X}(\mathbf{\Sigma^{e}})$ is
a Deligne-Mumford stack.
\end{prop}
\begin{pf}
From (\ref{bundle}),
$^{P}\mathcal{X}(\mathbf{\Sigma^{e}})=[(P\times_{(\mathbb{C}^{\times})^{m}}Z^{e})/G^{e}]$
is a quotient stack, where $G^{e}$ acts trivially on $P$. The
action of $G^{e}$ on $Z^{e}$ has finite, reduced stabilizers
because the stack $[Z^{e}/G^{e}]$ is a Deligne-Mumford stack, so
the action of $G^{e}$ on
$P\times_{(\mathbb{C}^{\times})^{m}}Z^{e}$ also has finite,
reduced stabilizers. From Corollary 2.2 of [9],
$^{P}\mathcal{X}(\mathbf{\Sigma^{e}})$ is a Deligne-Mumford stack.
\end{pf}

For an extended stacky fan $\mathbf{\Sigma^{e}}$, let $\sigma\in
\Sigma$ be a cone,  let $link(\tau):=\{\sigma: \sigma+\tau\in
\Sigma, \sigma\cap \tau=0\}$, and
$\widetilde{\rho}_{1},\ldots,\widetilde{\rho}_{l}$ be the rays in
$link(\sigma)$. Then
$\mathbf{\Sigma^{e}/\sigma}=(N(\sigma),\Sigma/\sigma,\beta^{e}(\sigma))$
is an extended stacky fan, where $\beta^{e}(\sigma):
\mathbb{Z}^{l+m-n}\longrightarrow N(\sigma)$ is given by the
images of $b_{1},\ldots,b_{l},b_{n+1},\ldots,b_{m}$ under
$N\longrightarrow N(\sigma)$. From the construction of extended
toric Deligne-Mumford stack, we have
$\mathcal{X}(\mathbf{\Sigma^{e}/\sigma}):=[Z^{e}(\sigma)/G^{e}(\sigma)]$,
where
$Z^{e}(\sigma)=(\mathbb{A}^{l}-\mathbb{V}(J_{\Sigma/\sigma}))\times
(\mathbb{C}^{\times})^{m-n}=Z(\sigma)\times
(\mathbb{C}^{\times})^{m-n}$,
$G^{e}(\sigma)=Hom_{\mathbb{Z}}(DG(\beta^{e}(\sigma)),\mathbb{C}^{\times})$.
We have an action of $(\mathbb{C}^{\times})^{m}$ on
$Z^{e}(\sigma)$ induced by the natural action of
$(\mathbb{C}^{\times})^{l+m-n}$ on $Z^{e}(\sigma)$ and the
projection $(\mathbb{C}^{\times})^{m}\longrightarrow
(\mathbb{C}^{\times})^{l+m-n}$. We let
\begin{eqnarray}
^{P}\mathcal{X}(\mathbf{\Sigma^{e}/\sigma})&=&[(P\times_{(\mathbb{C}^{\times})^{m}}(\mathbb{C}^{\times})^{l+m-n}
\times_{(\mathbb{C}^{\times})^{l+m-n}}Z^{e}(\sigma))/G^{e}(\sigma)] \nonumber \\
&=&[(P\times_{(\mathbb{C}^{\times})^{m}}Z^{e}(\sigma))/G^{e}(\sigma)]
\nonumber
\end{eqnarray}
be the quotient stack. Then we have:

\begin{prop}
Let $\sigma$ be a cone in the extended stacky fan
$\mathbf{\Sigma^{e}}$, then
$^{P}\mathcal{X}(\mathbf{\Sigma^{e}/\sigma})$ defines a closed
substack of $^{P}\mathcal{X}(\mathbf{\Sigma^{e}})$.
\end{prop}
\begin{pf}
Let $\mathcal{X}(\mathbf{\Sigma^{e}})=[Z^{e}/G^{e}]$,  if $\sigma$
is a cone, let $W^{e}(\sigma)$ be the closed subvariety of $Z^{e}$
defined by $J(\sigma):=<z_{i}:\rho_{i}\subseteq \sigma>$ in
$\mathbb{C}[z_{1},\ldots,z_{n},z_{n+1}^{\pm 1},\ldots,z_{m}^{\pm
1}]$, then we see that
$W^{e}(\sigma)=W(\sigma)\times(\mathbb{C}^{\times})^{m-n}$, where
$W(\sigma)$ is the closed subvariety of $Z$ defined by
$J(\sigma):=<z_{i}:\rho_{i}\subseteq \sigma>$ in
$\mathbb{C}[z_{1},\ldots,z_{n}]$.  From [4], there is a map
$\varphi_{0}: W(\sigma)\longrightarrow Z(\sigma)$ which is
$(\mathbb{C}^{\times})^{n}$-equivariant, we define the map
$W^{e}(\sigma)\longrightarrow Z^{e}(\sigma)$ by $\varphi_{0}\times
1$. Twist it by the bundle $P$, we have a map $\varphi_{0}:
P\times_{(\mathbb{C}^{\times})^{m}}W^{e}(\sigma)\longrightarrow
P\times_{(\mathbb{C}^{\times})^{m}}Z^{e}(\sigma)$. From the
following diagram:
\[
\begin{CD}
0 @ >>>\mathbb{Z}^{n-l}@ >>> \mathbb{Z}^{m}@ >>>
\mathbb{Z}^{l+m-n} @
>>> 0\\
&& @VV{}V@VV{\overline{\beta}^{e}}V@VV{\overline{\beta}^{e}(\sigma)}V \\
0@ >>>N_{\sigma} @ >{}>>N@ >>> N(\sigma) @>>> 0
\end{CD}
\]
Applying Gale dual and $Hom$ functor we get the commutative
diagram: \begin{equation}\Label{close} \vcenter{\xymatrix{
G^{e}\dto^{\alpha^{e}}\rto^{\varphi_{1}} &G^{e}(\sigma)\dto^{\alpha^{e}(\sigma)}\\
(\mathbb{C}^{\times})^{m}\rto^{} &(\mathbb{C}^{\times})^{l+m-n}}}
\end{equation}
So we get a map of groupoids: $\varphi_{0}\times\varphi_{1}:
P\times_{(\mathbb{C}^{\times})^{m}}W^{e}(\sigma)\times
G^{e}\longrightarrow
P\times_{(\mathbb{C}^{\times})^{m}}Z^{e}(\sigma)\times
G^{e}(\sigma)$ which is Morita equivalent.  So we have an
isomorphism of stacks
$[(P\times_{(\mathbb{C}^{\times})^{m}}W(\sigma))/G^{e}]\cong
[(P\times_{(\mathbb{C}^{\times})^{m}}Z^{e}(\sigma))/G^{e}(\sigma)]$.
Since $W^{e}(\sigma)$ is a subvariety of $Z^{e}$, and
$P\times_{(\mathbb{C}^{\times})^{m}}W^{e}(\sigma)$ is a subvariety
of $P\times_{(\mathbb{C}^{\times})^{m}}Z^{e}$, so
$[(P\times_{(\mathbb{C}^{\times})^{m}}W^{e}(\sigma))/G^{e}]$ is a
substack of
$[(P\times_{(\mathbb{C}^{\times})^{m}}Z^{e})/G^{e}]=~^{P}\mathcal{X}(\mathbf{\Sigma^{e}})$.
So $^{P}\mathcal{X}(\mathbf{\Sigma^{e}/\sigma})$ is a closed
substack of $^{P}\mathcal{X}(\mathbf{\Sigma^{e}})$.
\end{pf}

\begin{rmk}
From [4], $W(\sigma)=Z^{<g_{1},\cdots,g_{r}>}$ for some group
elements in $G$. From Proposition 2.3, the extended
Deligne-Mumford stack $[Z^{e}(\sigma)/G^{e}(\sigma)]$ is
isomorphic to the stack $[Z(\sigma)/G(\sigma)]$. Let
$g_{1},\cdots,g_{r}$ still represent the elements in $G^{e}$
through the map $\varphi_{1}$ in (\ref{second}). Then
$W^{e}(\sigma)=(Z^{e})^{<g_{1},\cdots,g_{r}>}$.
\end{rmk}

\begin{prop}
Let $^{P}\mathcal{X}(\mathbf{\Sigma^{e}})\longrightarrow B$ be a
toric stack bundle over a smooth variety $B$ with fibre
$\mathcal{X}(\mathbf{\Sigma^{e}})$ the  extended toric
Deligne-Mumford stack associated to the extended stacky fan
$\mathbf{\Sigma^{e}}$, then the $r$-th inertia stack of this toric
stack bundle is
$$\mathcal{I}_{r}\left(^{P}\mathcal{X}(\mathbf{\Sigma^{e}})\right)=\coprod_{(v_{1},\cdots,v_{r})\in Box(\mathbf{\Sigma^{e}})^{r}}
~^{P}\mathcal{X}(\mathbf{\Sigma^{e}/\sigma}(\overline{v}_{1},\cdots,\overline{v}_{r}))$$
\end{prop}
\begin{pf}
From (\ref{bundle}),
$^{P}\mathcal{X}(\mathbf{\Sigma^{e}})=[(P\times_{(\mathbb{C}^{\times})^{m}}Z^{e})/G^{e}]$
is a quotient stack. Because $G^{e}$ is an abelian group and the
the action has finite, reduced stabilizers, we have the $r$-th
inertia stack:
$$\mathcal{I}_{r}\left(^{P}\mathcal{X}(\mathbf{\Sigma^{e}})\right)=\left[\left(\coprod_{(g_{1},\cdots,g_{r})\in (G^{e})^{r}}
(P\times_{(\mathbb{C}^{\times})^{m}}Z^{e})^{H}\right)\diagup
G^{e}\right]$$ where $H$ is the subgroup in $G^{e}$ generated by
the elements $g_{1},\cdots,g_{r}$. From Lemma 4.6 in [4], there is
a map from $Box(\mathbf{\Sigma^{e}})$ to $G$, from the map
$\varphi_{1}$ in (\ref{second}), we have a map $\rho:
Box(\mathbf{\Sigma^{e}})\longrightarrow G^{e}$ such that
$\rho(v)=g(v)$. For a r-tuple $(v_{1},\cdots,v_{r})$ in the
$Box(\mathbf{\Sigma^{e}})$, from Proposition 3.5 and the above
Remark, we have:
$^{P}\mathcal{X}(\mathbf{\Sigma^{e}/\sigma}(\overline{v}_{1},\cdots,\overline{v}_{r}))\cong
[P\times_{(\mathbb{C}^{\times})^{m}}(Z^{e})^{H}/G^{e}]$. Taking
the disjoint union over all r-tuples in $Box(\mathbf{\Sigma^{e}})$
we get a map:
$$\psi: \coprod_{(v_{1},\cdots,v_{r})\in Box(\mathbf{\Sigma^{e}})^{r}}
~^{P}\mathcal{X}(\mathbf{\Sigma^{e}/\sigma}(\overline{v}_{1},\cdots,\overline{v}_{r}))
\longrightarrow
\mathcal{I}_{r}\left(^{P}\mathcal{X}(\mathbf{\Sigma^{e}})\right)$$
The toric stack bundle $^{P}\mathcal{X}(\mathbf{\Sigma^{e}})$
locally like a smooth variety times the extended toric
Deligne-Mumford stack $\mathcal{X}(\mathbf{\Sigma^{e}})$. From
[4], the map $\psi$ is an isomorphism locally in the Zariski
topology of the base $B$, so $\psi$ is an isomorphism  globally.
We complete the proof of the Proposition.
\end{pf}
\begin{rmk}
For any pair $(v_{1},v_{2})\in Box(\mathbf{\Sigma^{e}})^{2}$,
there exists a unique element $v_{3}\in Box(\mathbf{\Sigma^{e}})$
such that $v_{1}+v_{2}+v_{3}\equiv 0$. This means that in the
local group $N/N_{\sigma(\overline{v}_{1},\overline{v}_{2})}$, the
corresponding group elements $g_{1},g_{2},g_{3}$ satisfy
$g_{1}g_{2}g_{3}=1$. So this implies that
$\sigma(\overline{v}_{1},\overline{v}_{2},\overline{v}_{3})=\sigma(\overline{v}_{1},\overline{v}_{2})$.
In fact, the Proposition determines all the 3-twisted sectors of
the toric stack  bundle $^{P}\mathcal{X}(\mathbf{\Sigma^{e}})$.
See also in [17],[11] for the case of toric varieties.
\end{rmk}

%%% ----------------------------------------------------------------------
\section{The Orbifold Cohomology Ring.}
In this section we describe the ring structure of the orbifold
cohomology space of the toric stack bundles.
%%% ----------------------------------------------------------------------
\subsection{The Module Structure on $\bf{A^{*}_{orb}(^{P}\mathcal{X}(\mathbf{\Sigma^{e}})
)}$.}

Let $^{P}\mathcal{X}(\mathbf{\Sigma^{e}})\longrightarrow B$ be a
toric stack bundle. Let $E\longrightarrow B$ be the associated
$T$-bundle over $B$ induced from $P$ from the exact sequence
(\ref{stack}). Let $M:=N^{*}$ be the dual of $N$, and let
$\theta\in M$,  define $\xi_{\theta}\longrightarrow B$ to be the
line bundle coming from $E\longrightarrow B$ by "extending" the
structure group via $\chi^{\theta}: T\longrightarrow
\mathbb{C}^{\times}$. We give several definitions:
\begin{defn}
Let $A^{*}(B)$ denote the Chow ring over $\mathbb{Q}$ of the
smooth variety $B$. Define  the deformed ring
$A^{*}(B)[N]^{\mathbf{\Sigma^{e}}}$ as follows:
$A^{*}(B)[N]^{\mathbf{\Sigma^{e}}}=A^{*}(B)\otimes
_{\mathbb{Q}}\mathbb{Q}[N]^{\mathbf{\Sigma^{e}}}$,
$\mathbb{Q}[N]^{\mathbf{\Sigma^{e}}}=\bigoplus_{c\in N}y^{c}$,
where $y$ is  a formal variable. Multiplication is given by
(\ref{product}).
\end{defn}

The deformed ring $A^{*}(B)[N]^{\mathbf{\Sigma^{e}}}$ has a
$\mathbb{Q}$-grading defined by: if
$\overline{c}=\sum_{\rho_{i}\subseteq
\sigma(\overline{c})}a_{i}\overline{b}_{i}$, $deg(y^{c})=\sum
a_{i}\in \mathbb{Q}$. If $\gamma\in A^{*}(B)$, then
$deg(\gamma\cdot y^{c})=deg(\gamma)+deg(y^{c})$.  Let
$\mathcal{I}(\mathbf{\Sigma^{e}})$ be the ideal in (\ref{ideal}).

\begin{defn}
Let $\mathbf{\Sigma^{e}}=(N,\Sigma,\beta^{e})$ be an extended
stacky fan in $N_{\mathbb{Q}}$.  Define  ring
$S_{\mathbf{\Sigma^{e}}}:=A^{*}(B)[x_{1},\ldots,x_{n}]/I_{\mathbf{\Sigma^{e}}}$,
where the ideal $I_{\mathbf{\Sigma^{e}}}$ is generated by the
square-free monomials $x_{i_{1}}\cdots x_{i_{s}}$ with
$\rho_{i_{1}}+\cdots +\rho_{i_{s}}\notin \Sigma$.
\end{defn}

Note that $S_{\mathbf{\Sigma^{e}}}$ is a subring of
$A^{*}(B)[N]^{\mathbf{\Sigma^{e}}}$ given by the map
$x_{i}\longmapsto y^{b_{i}}$ for $1\leq i\leq n$.  Let
$\{\rho_{1},\ldots,\rho_{n}\}$ be the rays of
$\mathbf{\Sigma^{e}}$, then each $\rho_{i}$ corresponds to a line
bundle $L_{i}$ over the extended  toric Deligne-Mumford stack
$\mathcal{X}(\mathbf{\Sigma^{e}})$. This line bundle can be
defined as follow. The line bundle $L_{i}$ on the toric
Deligne-Mumford  stack $\mathcal{X}(\mathbf{\Sigma})$ is given by
the trivial line bundle $\mathbb{C}\times Z$ over $Z$ with the $G$
action on $\mathbb{C}$ given by the $i$-th component $\alpha_{i}$
of $\alpha: G\longrightarrow (\mathbb{C}^{\times})^{n}$ in
(\ref{stack}) when $\mathbf{\Sigma^{e}}=\mathbf{\Sigma}$. From
(\ref{second}), we have:
\begin{equation}\Label{pullback}\vcenter{\xymatrix{
G\dto^{\alpha}\rto^{\varphi_{1}} &G^{e}\dto^{\alpha^{e}}\\
(\mathbb{C}^{\times})^{n}\rto^{i} &(\mathbb{C}^{\times})^{m} }}
\end{equation}

\begin{defn}
For each  $\rho_{i}$, define the  line bundle $L_{i}$ over
$\mathcal{X}(\mathbf{\Sigma^{e}})$ to be the quotient of the
trivial line bundle $Z^{e}\times \mathbb{C}$ over $Z^{e}$ under
the action of $G^{e}$ on $\mathbb{C}$ through one component of
$\alpha^{e}$ such that the pullback component in $\alpha$ through
(\ref{pullback}) is $\alpha_{i}$. Twist it by the principal
$(\mathbb{C}^{\times})^{m}$-bundle $P$, we get the line bundle
$\mathcal{L}_{i}$ over the toric stack bundle
$^{P}\mathcal{X}(\mathbf{\Sigma^{e}})$.
\end{defn}
First we describe the ordinary Chow ring of the toric stack
bundle:

\begin{lem}
Let $^{P}\mathcal{X}(\mathbf{\Sigma^{e}})\longrightarrow B$ be a
toric stack bundle over a smooth variety $B$ with  fibre
$\mathcal{X}(\mathbf{\Sigma^{e}})$ the extended toric
Deligne-Mumford stack associated to the extended stacky fan
$\mathbf{\Sigma^{e}}$, then there is an isomorphism of
$\mathbb{Q}$-graded rings:
$$\frac{S_{\mathbf{\Sigma^{e}}}}{\mathcal{I}(\mathbf{\Sigma^{e}})}\cong A^{*}\left(^{P}\mathcal{X}(\mathbf{\Sigma^{e}})\right)$$
given by $x_{i}\longmapsto c_{1}(\mathcal{L}_{i})$.
\end{lem}
\begin{pf}
From Corollary 2.4, let $X(\Sigma)$ be the coarse moduli space of
the extended toric Deligne-Mumford stack
$\mathcal{X}(\mathbf{\Sigma^{e}})$. Let $E\longrightarrow B$ be
the principal $T$-bundle induced from the
$(\mathbb{C}^{\times})^{m}$-bundle $P$, then from Proposition 3.3,
$^{E}X(\Sigma)$ is the coarse moduli space of the toric stack
bundle $^{P}\mathcal{X}(\mathbf{\Sigma^{e}})$. Let $a_{i}$ be the
first lattice vector in the ray generated by $b_{i}$, then
$\overline{b}_{i}=l_{i}a_{i}$ for some positive integer $l_{i}$.
The ideal $\mathcal{I}(\mathbf{\Sigma^{e}})$ in (\ref{ideal}) also
define an ideal in $S_{\mathbf{\Sigma^{e}}}$. From [19], we have
$$A^{*}(^{E}X(\Sigma))\cong
\frac{S_{\mathbf{\Sigma^{e}}}}{\mathcal{I}(\mathbf{\Sigma^{e}})}$$
which is given by $x_{i}\longmapsto E(V(\rho_{i}))$, where
$E(V(\rho_{i}))$ is the associated bundle over $B$ corresponding
to the $T$-invariant divisor $V(\rho_{i})$. From [2],[20], the
Chow ring of the stack $^{P}\mathcal{X}(\mathbf{\Sigma^{e}})$ is
isomorphic to the Chow ring of its coarse moduli space
$^{E}X(\Sigma)$ given by $c_{1}(\mathcal{L}_{i})\longmapsto
l_{i}^{-1}\cdot E(V(\rho_{i}))$, and
$c_{1}(\xi_{\theta})+\sum_{i=1}^{n}\theta(a_{i})l_{i}y^{b_{i}}=c_{1}(\xi_{\theta})+\sum_{i=1}^{n}\theta(b_{i})y^{b_{i}}$,
so we prove the Lemma.
\end{pf}

Now we talk about the module structure on
$A_{orb}^{*}\left(^{P}\mathcal{X}(\mathbf{\Sigma^{e}})\right)$.
Because $\Sigma$ is a simplicial fan, we have:
\begin{lem}
For any $c\in N$, let $\sigma$ be the minimal cone in $\Sigma$
containing $\overline{c}$, then there exists a unique expression
$$c=v+\sum_{\rho_{i}\subset\sigma}m_{i}b_{i}$$
where $m_{i}\in \mathbb{Z}_{\geq 0}$, and $v\in Box(\sigma)$.
~$\square$
\end{lem}

\begin{lem}
If $\tau$ is a cone in the complete simplicial fan $\Sigma$,
$\{\rho_{1},\ldots,\rho_{s}\}\subset link(\tau)$, suppose
$\rho_{1},\ldots,\rho_{s}$ are contained in a cone $\sigma\subset
\Sigma$. Then $\sigma\cup\tau$ is contained in a cone of $\Sigma$.
\end{lem}
\begin{pf}
Using the following result: If $\rho_{1},\ldots,\rho_{s}$ are rays
in the complete simplicial fan $\Sigma$,  if for any $i,j$,
$\rho_{i},\rho_{j}$ generate a cone, then
$\rho_{1},\cdots,\rho_{s}$ generate a cone, see [10],[15]. The
Lemma is proved.
\end{pf}

\begin{prop}
Let $^{P}\mathcal{X}(\mathbf{\Sigma^{e}})\longrightarrow B$ be a
toric stack bundle over a smooth variety $B$ with fibre
$\mathcal{X}(\mathbf{\Sigma^{e}})$ the extended toric
Deligne-Mumford stack associated to the extended stacky fan
$\mathbf{\Sigma^{e}}$, then we have an isomorphism of
$A^{*}(^{P}\mathcal{X}(\mathbf{\Sigma^{e}}))$-modules:
$$\bigoplus_{v\in Box(\mathbf{\Sigma^{e}})}A^{*}\left(^{P}\mathcal{X}(\mathbf{\Sigma^{e}/\sigma}(\overline{v}))\right)[deg(y^{v})]\cong
\frac{A^{*}(B)[N]^{\mathbf{\mathbf{\Sigma^{e}}}}}{\mathcal{I}(\mathbf{\Sigma^{e}})}$$
\end{prop}
\begin{pf}
From the definition of $A^{*}(B)[N]^{\mathbf{\Sigma^{e}}}$ and
Lemma 4.5, we see  that
$A^{*}(B)[N]^{\mathbf{\Sigma^{e}}}=\bigoplus_{v\in
Box(\mathbf{\Sigma^{e}})}y^{v}\cdot S_{\mathbf{\Sigma^{e}}}$.
Since $\mathcal{I}(\mathbf{\Sigma^{e}})$ is the ideal in
$A^{*}(B)[N]^{\mathbf{\Sigma^{e}}}$ defined in (\ref{ideal}). Then
$\bigoplus_{v\in Box(\mathbf{\Sigma^{e}})}y^{v}\cdot
\mathcal{I}(\mathbf{\Sigma^{e}})$ is the ideal
$\mathcal{I}(\mathbf{\Sigma^{e}})$ in $\bigoplus_{v\in
Box(\mathbf{\Sigma^{e}})}y^{v}\cdot
S_{\mathbf{\Sigma^{e}}}=A^{*}(B)[N]^{\mathbf{\Sigma^{e}}}$.  So we
obtain the isomorphism of
$A^{*}(^{P}\mathcal{X}(\mathbf{\Sigma^{e}}))$-modules:
\begin{equation}\Label{module}
\frac{A^{*}(B)[N]^{\mathbf{\mathbf{\Sigma^{e}}}}}{\mathcal{I}(\mathbf{\Sigma^{e}})}\cong
\bigoplus_{v\in Box(\mathbf{\Sigma^{e}})}\frac{y^{v}\cdot
S_{\mathbf{\Sigma^{e}}}}{y^{v}\cdot
\mathcal{I}(\mathbf{\Sigma^{e}})}
\end{equation}

For any  $v\in Box(\mathbf{\Sigma^{e}})$,  let
$\sigma(\overline{v})$ be the minimal cone in $\Sigma$ containing
$\overline{v}$.   Let $\rho_{1},\ldots, \rho_{l}\in
link(\sigma(\overline{v}))$, and $\widetilde{\rho}_{i}$ be the
image of $\rho_{i}$ under the natural map $N\longrightarrow
N(\sigma(\overline{v}))=N/N_{\sigma(\overline{v})}$. Then
$S_{\mathbf{\Sigma^{e}}/\sigma(\overline{v})}\subset
A^{*}(B)[N(\sigma(\overline{v}))]^{\mathbf{\Sigma^{e}}/\sigma(\overline{v})}$
is the subring given by: $\widetilde{x}_{i}\longmapsto
y^{\widetilde{b}_{i}}$, for $\rho_{i}\in
link(\sigma(\overline{v}))$. Consider the morphism: $i:
A^{*}(B)[\widetilde{x}_{1},\ldots,\widetilde{x}_{l}]\longrightarrow
A^{*}(B)[x_{1},\ldots,x_{n}]$ given by
$\widetilde{x}_{i}\longrightarrow x_{i}$. From Lemma 4.6, it is
easy to check that the ideal
$I_{\mathbf{\Sigma^{e}}/\sigma(\overline{v})}$  goes to the ideal
$I_{\mathbf{\Sigma^{e}}}$,  so we have a  morphism
$S_{\mathbf{\Sigma^{e}}/\sigma(\overline{v})}\longrightarrow
S_{\mathbf{\Sigma^{e}}}$. Since $S_{\mathbf{\Sigma^{e}}}$ is a
subring of $A^{*}(B)[N]^{\mathbf{\Sigma^{e}}}$ given by
$x_{i}\longmapsto y^{b_{i}}$, we use the notations $y^{b_{i}}$.
Let $\widetilde{\Psi}_{v}:
S_{\mathbf{\Sigma^{e}}/\sigma(\overline{v})}[deg(y^{v})]\longrightarrow
y^{v}\cdot S_{\mathbf{\Sigma^{e}}}$ be the morphism given by:
$y^{\widetilde{b}_{i}}\longmapsto y^{v}\cdot y^{b_{i}}$. If
$\sum_{i=1}^{l}\widetilde{\theta}(\widetilde{b}_{i})y^{\widetilde{b}_{i}}+c_{1}(\xi_{\widetilde{\theta}})$
belongs to the ideal
$\mathcal{I}(\mathbf{\Sigma^{e}/\sigma}(\overline{v}))$, then
$$
\widetilde{\Psi}_{v}\left(\sum_{i=1}^{l}\widetilde{\theta}(\widetilde{b}_{i})y^{\widetilde{b}_{i}}+c_{1}(\xi_{\widetilde{\theta}})\right)
= y^{v}\cdot
\left(\sum_{i=1}^{l}\widetilde{\theta}(\widetilde{b}_{i})y^{\widetilde{b}_{i}}+c_{1}(\xi_{\widetilde{\theta}})\right)
=y^{v}\cdot
\left(\sum_{i=1}^{n}\theta(b_{i})y^{b_{i}}+c_{1}(\xi_{\theta})\right)
$$
where $\theta$ is determined by the diagram:
\begin{equation}\Label{map}
\vcenter{\xymatrix{
N\dto_{\pi}\drto^{\theta} & \\
N(\sigma(\overline{v}))\rto^{\widetilde{\theta}} & \ \mathbb{Z}}}
\end{equation}
So $\theta(b_{i})=\widetilde{\theta}(\widetilde{b}_{i})$. From the
definition  of the line bundle $\xi_{\theta}$, we have
$\xi_{\theta}\cong \xi_{\widetilde{\theta}}$. We obtain that
$\widetilde{\Psi}_{v}(\sum_{i=1}^{l}\widetilde{\theta}(\widetilde{b}_{i})y^{\widetilde{b}_{i}}+c_{1}(\xi_{\widetilde{\theta}}))\in
y^{v}\cdot \mathcal{I}(\mathbf{\Sigma^{e}})$. So
$\widetilde{\Psi}_{v}$ induce  a morphism $\Psi_{v}:
\frac{S_{\mathbf{\Sigma^{e}}/\sigma(\overline{v})}}{\mathcal{I}(\mathbf{\Sigma^{e}}/\sigma(\overline{v}))}[deg(y^{v})]
\longrightarrow \frac{y^{v}\cdot
S_{\mathbf{\Sigma^{e}}}}{y^{v}\cdot
\mathcal{I}(\mathbf{\Sigma^{e}})}
$
such that $\Psi_{v}([y^{\widetilde{b}_{i}}])=[y^{v}\cdot
y^{b_{i}}]$.

Conversely, for such $v\in Box(\mathbf{\Sigma^{e}})$ and
$\rho_{i}\subset \sigma(\overline{v})$, choose $\theta_{i}\in
Hom(N,\mathbb{Q})$ such that $\theta_{i}(b_{i})=1$ and
$\theta_{i}(b_{i^{'}})=0$ for $b_{i^{'}}\neq b_{i}\in
\sigma(\overline{v})$.  We  consider the following morphism $p:
A^{*}(B)[x_{1},\ldots,x_{n}]\longrightarrow
A^{*}(B)[\widetilde{x}_{1},\ldots,\widetilde{x}_{l}]$, where $p$
is given by:
$$x_{i}\longmapsto\begin{cases}\widetilde{x}_{i}&\text{if $\rho_{i}\subseteq link(\sigma(\overline{v}))$}\,,\\
-\sum_{j=1}^{l}\theta_{i}(b_{j})\widetilde{x}_{j}&\text{if
$\rho_{i}\subseteq \sigma(\overline{v})$}\,,\\
0&\text{if $\rho_{i}\nsubseteq \sigma(\overline{v})\cup
link(\sigma(\overline{v}))$}\,.\end{cases}$$  For any
$x_{i_{1}}\cdots x_{i_{s}}$ in $I_{\mathbf{\Sigma^{e}}}$, also
from Lemma 4.6 we prove that  $p(x_{i_{1}}\cdots x_{i_{s}})\in
I_{\mathbf{\Sigma^{e}}/\sigma(\overline{v})}$.  We also use the
notations $y^{b_{i}}$ to replace $x_{i}$, then $p$ induces a
surjective map: $S_{\mathbf{\Sigma^{e}}}\longrightarrow
S_{\mathbf{\Sigma^{e}}/\sigma(\overline{v})}$ and a surjective
map: $\widetilde{\Phi}_{v}: y^{v}\cdot
S_{\mathbf{\Sigma^{e}}}\longrightarrow
S_{\mathbf{\Sigma^{e}}/\sigma(\overline{v})}[deg(y^{v})]$. Let
$y^{v}\cdot
\left(\sum_{i=1}^{n}\theta(b_{i})y^{b_{i}}+c_{1}(\xi_{\theta})\right)$
belong to the ideal $y^{v}\cdot \mathcal{I}(\mathbf{\Sigma^{e}})$.
For $\theta\in M$, we have $\theta=\theta_{v}+\theta_{v}^{'}$,
where $\theta_{v}\in N(\sigma(\overline{v}))^{*}=M\cap
\sigma(\overline{v})^{\perp}$ and  $\theta_{v}^{'}$ belongs to the
orthogonal complement of the subspace
$\sigma(\overline{v})^{\perp}$ in $M$.  From (\ref{map}), we have:
\begin{eqnarray}
&~&\widetilde{\Phi}_{v}\left(y^{v}\cdot
\left(\sum_{i=1}^{n}\theta(b_{i})y^{b_{i}}+c_{1}(\xi_{\theta})\right)\right)
\nonumber \\
&=&
\sum_{i=1}^{l}\theta_{v}(\widetilde{b}_{i})y^{\widetilde{b}_{i}}+c_{1}(\xi_{\theta_{v}})+
\sum_{\rho_{i}\subset
\sigma(\overline{v})}\theta_{v}^{'}(b_{i})\left(-\sum_{j=1}^{l}\theta_{i}(b_{j})y^{\widetilde{b}_{j}}\right)+
c_{1}(\xi_{\theta_{v}^{'}})+
\sum_{i=1}^{l}\theta_{v}^{'}(b_{i})y^{\widetilde{b}_{i}} \nonumber
\end{eqnarray}
Note that
$\left(\sum_{i=1}^{l}\theta_{v}(\widetilde{b}_{i})y^{\widetilde{b}_{i}}+c_{1}(\xi_{\theta_{v}})\right)\in
\mathcal{I}(\mathbf{\Sigma^{e}}/\sigma(\overline{v}))$. From the
definition of $\xi_{\widetilde{\theta}}$  over
$\mathcal{X}(\mathbf{\Sigma^{e}}/\sigma(\overline{v}))$,
$\xi_{\theta_{v}^{'}}=0$. Now let
$\theta_{v}^{'}=\sum_{\rho_{i}\subset
\sigma(\overline{v})}a_{i}\theta_{i}$, where $a_{i}\in
\mathbb{Q}$,  then $\sum_{\rho_{i}\subset
\sigma(\overline{v})}\theta_{v}^{'}(b_{i})=\sum_{\rho_{i}\subset
\sigma(\overline{v})}a_{i}\theta_{i}(b_{i})$. We have:
$\sum_{\rho_{i}\subset
\sigma(\overline{v})}a_{i}\theta_{i}(b_{i})\left(-\sum_{j=1}^{l}\theta_{i}(b_{j})y^{\widetilde{b}_{j}}\right)
+\sum_{\rho_{i}\subset
\sigma(\overline{v})}\sum_{j=1}^{l}a_{i}\theta_{i}(b_{j})y^{\widetilde{b}_{j}}=0$,
so we have $\widetilde{\Phi}_{v}\left(y^{v}\cdot
\left(\sum_{i=1}^{n}\theta(b_{i})y^{b_{i}}+c_{1}(\xi_{\theta})\right)\right)\in
\mathcal{I}(\mathbf{\Sigma^{e}}/\sigma(\overline{v}))$. So
$\widetilde{\Phi}_{v}$ induces a morphism $\Phi: \frac{y^{v}\cdot
S_{\mathbf{\Sigma^{e}}}}{y^{v}\cdot
\mathcal{I}(\mathbf{\Sigma^{e}})}\longrightarrow
\frac{S_{\mathbf{\Sigma^{e}/\sigma}(\overline{v})}}{\mathcal{I}(\mathbf{\Sigma^{e}}/\sigma(\overline{v}))}[deg(y^{v})]$.
Note that $\Phi_{v}\Psi_{v}=1$ is easy to check. For any
$[y^{v}\cdot y^{b_{i}}]\in \frac{y^{v}\cdot
S_{\mathbf{\Sigma^{e}}}}{y^{v}\cdot
\mathcal{I}(\mathbf{\Sigma^{e}})}$, since $y^{v}\cdot
\left(-\sum_{j=1}^{l}\theta_{i}(b_{j})y^{b_{j}}+\sum_{j=1}^{n}\theta_{i}(b_{j})y^{b_{j}}\right)=y^{v}\cdot
y^{b_{i}}$,  we have
$[y^{v}\cdot(-\sum_{j=1}^{l}\theta_{i}(b_{j})y^{b_{j}})]=[y^{v}\cdot
y^{b_{i}}]$, so we check that $\Psi_{v}\Phi_{v}=1$. So $\Phi_{v}$
is an isomorphism.  From Lemma 4.4, for any $v\in
Box(\mathbf{\Sigma^{e}})$, we have an isomorphism of Chow rings:
$\frac{S_{\mathbf{\Sigma^{e}}/\sigma(\overline{v})}}{\mathcal{I}(\mathbf{\Sigma^{e}}/\sigma(\overline{v}))}\cong
A^{*}(^{P}\mathcal{X}(\mathbf{\Sigma^{e}}/\sigma(\overline{v})))$.
Taking into account all the $v$ in $Box(\mathbf{\Sigma^{e}})$ and
(\ref{module}) we have the isomorphism: $\bigoplus_{v\in
Box(\mathbf{\Sigma^{e}})}A^{*}\left(^{P}\mathcal{X}(\mathbf{\Sigma^{e}/\sigma}(\overline{v}))\right)[deg(y^{v})]\cong
\frac{A^{*}(B)[N]^{\mathbf{\mathbf{\Sigma^{e}}}}}{\mathcal{I}(\mathbf{\Sigma^{e}})}$.
Note that both sides of (\ref{module}) are
$S_{\mathbf{\Sigma^{e}}}/\mathcal{I}(\mathbf{\Sigma^{e}})=A^{*}(^{P}\mathcal{X}(\mathbf{\Sigma^{e}}))$-modules,
we complete the proof.
\end{pf}
\begin{rmk}
In Proposition 5.2 of [4], the authors give a proof of the
 Proposition for toric Deligne-Mumford stacks. We give a more
explicit proof of this isomorphism for the toric stack bundle in
this  Proposition.
\end{rmk}

%%% ----------------------------------------------------------------------
\subsection{The Orbifold Cup Product.}
In this section we consider the orbifold cup product on
$A_{orb}^{*}\left(^{P}\mathcal{X}(\mathbf{\Sigma^{e}})\right)$.
First we determine the 3-twisted sectors of
$^{P}\mathcal{X}(\mathbf{\Sigma^{e}})$. From the orbifold
Gromov-Witten theory, the 3-twisted sectors of
$^{P}\mathcal{X}(\mathbf{\Sigma^{e}})$ are the components of the
double inertia stack
$\mathcal{I}_{2}(^{P}\mathcal{X}(\mathbf{\Sigma^{e}}))$ of
$^{P}\mathcal{X}(\mathbf{\Sigma^{e}})$, see [6]. So from the
remark after Proposition 3.6, we have all the 3-twisted sectors of
$^{P}\mathcal{X}(\mathbf{\Sigma^{e}})$:
\begin{equation}\Label{3-sector}
\coprod_{(g_{1},g_{2},g_{3})\in Box(\mathbf{\Sigma^{e}})^{3},
g_{1}g_{2}g_{3}=1}
~^{P}\mathcal{X}(\mathbf{\Sigma^{e}/\sigma}(\overline{g}_{1},\overline{g}_{2},\overline{g}_{3}))
\end{equation}
where $\sigma(\overline{g}_{1},\overline{g}_{2},\overline{g}_{3})$
is the minimal cone in $\Sigma$ containing
$\overline{g}_{1},\overline{g}_{2},\overline{g}_{3}$. For any
3-twisted sector
$^{P}\mathcal{X}(\mathbf{\Sigma^{e}})_{(g_{1},g_{2},g_{3})}=
~^{P}\mathcal{X}(\mathbf{\Sigma^{e}/\sigma}(\overline{g}_{1},\overline{g}_{2},\overline{g}_{3}))$,
we have an inclusion
$e:~^{P}\mathcal{X}(\mathbf{\Sigma^{e}/\sigma}(\overline{g}_{1},\overline{g}_{2},\overline{g}_{3}))
\longrightarrow ~^{P}\mathcal{X}(\mathbf{\Sigma^{e}})$  because
$~^{P}\mathcal{X}(\mathbf{\Sigma^{e}/\sigma}(\overline{g}_{1},\overline{g}_{2},\overline{g}_{3}))$
is a substack of $~^{P}\mathcal{X}(\mathbf{\Sigma^{e}})$. Let $H$
be the subgroup generated by $g_{1},g_{2},g_{3}$, then the genus
zero, degree zero orbifold stable map to
$^{P}\mathcal{X}(\mathbf{\Sigma^{e}})$ determines a Galois
covering $\pi: C\longrightarrow \mathbb{P}^{1}$ branching over
three marked points $0,1,\infty$ such that the transformation
group of this covering is $H$. We have the definition:
\begin{defn}([5])
The obstruction bundle $O_{(g_{1},g_{2},g_{3})}$ over
$~^{P}\mathcal{X}(\mathbf{\Sigma^{e}/\sigma}(\overline{g}_{1},\overline{g}_{2},\overline{g}_{3}))$
is defined by the $H$-invariant:
$$\left(e^{*}T\left(^{P}\mathcal{X}(\mathbf{\Sigma^{e}})\right)\otimes
H^{1}(C,\mathcal{O}_{C})\right)^{H}$$
\end{defn}

\begin{prop}
Let $^{P}\mathcal{X}(\mathbf{\Sigma^{e}})_{(g_{1},g_{2},g_{3})}=
~^{P}\mathcal{X}(\mathbf{\Sigma^{e}/\sigma}(\overline{g}_{1},\overline{g}_{2},\overline{g}_{3}))$
be a 3-twisted sector of the stack
$^{P}\mathcal{X}(\mathbf{\Sigma^{e}})$, let
$g_{1}+g_{2}+g_{3}=\sum_{\rho_{i}\subset
\sigma(\overline{g}_{1},\overline{g}_{2},\overline{g}_{3})}a_{i}b_{i}$,
$a_{i}=1 , 2$, then the Euler class of the obstruction bundle
$O_{(g_{1},g_{2},g_{3})}$  over
$^{P}\mathcal{X}(\mathbf{\Sigma^{e}})_{(g_{1},g_{2},g_{3})}$ is:
$$\prod_{a_{i}=2}c_{1}(\mathcal{L}_{i})|_{^{P}\mathcal{X}(\mathbf{\Sigma^{e}/\sigma}(\overline{g}_{1},\overline{g}_{2},\overline{g}_{3}))}$$
where $\mathcal{L}_{i}$ is the line bundle over
$^{P}\mathcal{X}(\mathbf{\Sigma^{e}})$ in definition 4.3.
\end{prop}
\begin{pf}
Let $\mathcal{X}(\mathbf{\Sigma^{e}})$ be the extended toric
Deligne-Mumford stack corresponding to the extended  stacky fan
$\mathbf{\Sigma^{e}}$ in
$N_{\mathbb{Q}}=N\otimes_{\mathbb{Z}}\mathbb{Q}$. Let
$\sigma(\overline{g}_{1},\overline{g}_{2},\overline{g}_{3})$ be
the minimal cone in $\Sigma$ containing
$\overline{g}_{1},\overline{g}_{2},\overline{g}_{3}$.  From
Corollary 2.5 and (\ref{3-sector}) we have the 3-twisted sector
$\mathcal{X}(\mathbf{\Sigma^{e}})_{(g_{1},g_{2},g_{3})}=\mathcal{X}(\mathbf{\Sigma^{e}}/\sigma(\overline{g}_{1},\overline{g}_{2},\overline{g}_{3}))$
and $^{P}\mathcal{X}(\mathbf{\Sigma^{e}})_{(g_{1},g_{2},g_{3})}
=~^{P}\mathcal{X}(\mathbf{\Sigma^{e}/\sigma}(\overline{g}_{1},\overline{g}_{2},\overline{g}_{3}))$.
Since $e:
\mathcal{X}(\mathbf{\Sigma^{e}})_{(g_{1},g_{2},g_{3})}\longrightarrow
\mathcal{X}(\mathbf{\Sigma^{e}})$  is an inclusion, we have an
exact sequence:
$$0\longrightarrow T\mathcal{X}(\mathbf{\Sigma^{e}/\sigma}(\overline{g}_{1},\overline{g}_{2},\overline{g}_{3}))
\longrightarrow e^{*}T\mathcal{X}(\mathbf{\Sigma^{e}})
\longrightarrow
N(\mathcal{X}(\mathbf{\Sigma^{e}/\sigma}(\overline{g}_{1},\overline{g}_{2},\overline{g}_{3}))/\mathcal{X}(\mathbf{\Sigma^{e}}))\longrightarrow
0$$ where
$N(\mathcal{X}(\mathbf{\Sigma^{e}/\sigma}(\overline{g}_{1},\overline{g}_{2},\overline{g}_{3}))/\mathcal{X}(\mathbf{\Sigma^{e}}))$
is the normal bundle of
$\mathcal{X}(\mathbf{\Sigma^{e}/\sigma}(\overline{g}_{1},\overline{g}_{2},\overline{g}_{3}))$
in $\mathcal{X}(\mathbf{\Sigma^{e}})$.

Since $\mathcal{X}(\mathbf{\Sigma^{e}}):=[Z^{e}/G^{e}]$, the
tangent bundle
$T(\mathcal{X}(\mathbf{\Sigma^{e}}))=[T(Z^{e})/T(G^{e})]$ is a
quotient stack. Since $Z^{e}$ is an open subvariety of
$\mathbb{A}^{n}\times (\mathbb{C}^{\times})^{m-n}$,
$T(Z^{e})=\mathcal{O}_{Z^{e}}^{n}$. Now from the construction of
the line bundle $L_{k}$ over $\mathcal{X}(\mathbf{\Sigma^{e}})$,
we have a canonical map: $\bigoplus_{k=1}^{n}L_{k}\longrightarrow
T(\mathcal{X}(\mathbf{\Sigma^{e}}))$. Since we have a natural map
$T(\mathcal{X}(\mathbf{\Sigma^{e}}))\longrightarrow
N(\mathcal{X}(\mathbf{\Sigma^{e}/\sigma}(\overline{g}_{1},\overline{g}_{2},\overline{g}_{3}))/\mathcal{X}(\mathbf{\Sigma^{e}}))$,
we obtain a map of vector bundles over
$\mathcal{X}(\mathbf{\Sigma^{e}/\sigma}(\overline{g}_{1},\overline{g}_{2},\overline{g}_{3}))$:
$$\varphi: \bigoplus_{\rho_{k}\subset
\sigma(\overline{g}_{1},\overline{g}_{2},\overline{g}_{3})}L_{k}\longrightarrow
N(\mathcal{X}(\mathbf{\Sigma^{e}/\sigma}(\overline{g}_{1},\overline{g}_{2},\overline{g}_{3}))/\mathcal{X}(\mathbf{\Sigma^{e}}))$$
Then from the definition of the line bundle $\mathcal{L}_{k}$ over
$^{P}\mathcal{X}(\mathbf{\Sigma^{e}})$, we have the map:
$$\widetilde{\varphi}: \bigoplus_{\rho_{k}\subset
\sigma(\overline{g}_{1},\overline{g}_{2},\overline{g}_{3})}\mathcal{L}_{k}\longrightarrow
N(^{P}\mathcal{X}(\mathbf{\Sigma^{e}/\sigma}(\overline{g}_{1},\overline{g}_{2},\overline{g}_{3}))/^{P}\mathcal{X}(\mathbf{\Sigma^{e}}))$$
where
$N(^{P}\mathcal{X}(\mathbf{\Sigma^{e}/\sigma}(\overline{g}_{1},\overline{g}_{2},\overline{g}_{3}))/^{P}\mathcal{X}(\mathbf{\Sigma^{e}}))$
is the normal bundle of
$^{P}\mathcal{X}(\mathbf{\Sigma^{e}/\sigma}(\overline{g}_{1},\overline{g}_{2},\overline{g}_{3}))$
in $^{P}\mathcal{X}(\mathbf{\Sigma^{e}})$.  For any point map:
$$x: Spec~ \mathbb{C}\hookrightarrow \mathcal{X}(\mathbf{\Sigma^{e}/\sigma}(\overline{g}_{1},\overline{g}_{2},\overline{g}_{3}))\hookrightarrow
~^{P}\mathcal{X}(\mathbf{\Sigma^{e}/\sigma}(\overline{g}_{1},\overline{g}_{2},\overline{g}_{3}))$$
note that $x^{*}\widetilde{\varphi}$ is an isomorphism, so
$\widetilde{\varphi}$ is an isomorphism. We have the exact
sequence:
$$0\longrightarrow T\left(^{P}\mathcal{X}(\mathbf{\Sigma^{e}/\sigma}(\overline{g}_{1},\overline{g}_{2},\overline{g}_{3}))\right)\longrightarrow
e^{*}T\left(^{P}\mathcal{X}(\mathbf{\Sigma^{e}})\right)
\longrightarrow \bigoplus_{\rho_{k}\subset
\sigma(\overline{g}_{1},\overline{g}_{2},\overline{g}_{3})}\mathcal{L}_{k}\longrightarrow
0$$ Now using the result in the proof of Proposition 6.3 in [4],
we have
$$(dim_{\mathbb{C}}(\mathcal{L}_{k}\otimes H^{1}(C,\mathcal{O}_{C}))^{H}=0 ~\mbox{if}~~
a_{k}=1,~~~ (dim_{\mathbb{C}}(\mathcal{L}_{k}\otimes
H^{1}(C,\mathcal{O}_{C}))^{H}=1 ~\mbox{if}~~ a_{k}=2$$ So from the
definition 4.8, we have:
$$e(O_{(g_{1},g_{2},g_{3})})\cong
\prod_{a_{i}=2}c_{1}(\mathcal{L}_{i})|_{^{P}\mathcal{X}(\mathbf{\Sigma^{e}/\sigma}(\overline{g}_{1},\overline{g}_{2},\overline{g}_{3}))}$$
\end{pf}

\subsection*{Proof of Theorem 1.1:}
From the definition of the orbifold cohomology  in [5] and
Proposition 4.7,  we know that
$A^{*}_{orb}\left(^{P}\mathcal{X}(\mathbf{\Sigma^{e}})\right)=\bigoplus_{g\in
Box(\mathbf{\Sigma^{e}})}A^{*}\left(^{P}\mathcal{X}(\mathbf{\Sigma^{e}/\sigma}(\overline{g}))\right)[deg(y^{g})]$,
and from Proposition 4.7, we have an isomorphism between
$A^{*}(^{P}\mathcal{X}(\mathbf{\Sigma^{e}}))$-modules:
$$\bigoplus_{g\in
Box(\mathbf{\Sigma^{e}})}A^{*}\left(^{P}\mathcal{X}(\mathbf{\Sigma^{e}/\sigma}(\overline{g}))\right)[deg(y^{g})]\cong
A^{*}(B)[N]^{\mathbf{\Sigma^{e}}}/\mathcal{I}(\mathbf{\Sigma^{e}})$$
So we have an isomorphism of
$A^{*}(^{P}\mathcal{X}(\mathbf{\Sigma^{e}}))$-modules:
$A^{*}_{orb}\left(^{P}\mathcal{X}(\mathbf{\Sigma^{e}})\right)\cong\frac{A^{*}(B)[N]^{\mathbf{\Sigma^{e}}}}{\mathcal{I}(\mathbf{\Sigma^{e}})}$.
Next we show that the orbifold cup product defined in [5]
coincides with the product in ring
$A^{*}(B)[N]^{\mathbf{\Sigma^{e}}}/\mathcal{I}(\mathbf{\Sigma^{e}})$.
From the above isomorphisms, it suffices to consider the canonical
generators $y^{b_{i}}$, $y^{g}$ where $g\in
Box(\mathbf{\Sigma^{e}})$ and $\gamma\in A^{*}(B)$. Since
$b_{i}\in N$, the twisted sector determined by $b_{i}$ is  the
whole toric stack bundle $^{P}\mathcal{X}(\mathbf{\Sigma^{e}})$,
$y^{b_{i}}\cup_{orb}\gamma$ is the usual product
$y^{b_{i}}\cdot\gamma$ in the deformed ring because $y^{b_{i}}$
and $\gamma$ belong to the ordinary Chow ring of
$^{P}\mathcal{X}(\mathbf{\Sigma^{e}})$.

For $y^{g}\cup_{orb}y^{b_{i}}$ and $y^{g}\cup_{orb}\gamma$, $g\in
Box(\mathbf{\Sigma^{e}})$, so $g$ determine a twisted sector
$^{P}\mathcal{X}(\mathbf{\Sigma^{e}/\sigma}(\overline{g}))$. The
corresponding twisted sector to $b_{i}$ and $\gamma$ are the whole
toric stack bundle $^{P}\mathcal{X}(\mathbf{\Sigma^{e}})$. It is
easy to see that the 3-twisted sector corresponding to $(g,b_{i})$
and $(g,\gamma)$ are
$^{P}\mathcal{X}(\mathbf{\Sigma^{e}})_{(g,1,g^{-1})}\cong
~^{P}\mathcal{X}(\mathbf{\Sigma^{e}/\sigma}(\overline{g}))$, where
$g^{-1}$ is the inverse of $g$ in the local group. From the
dimension formula in [5], the obstruction bundle over
$^{P}\mathcal{X}(\mathbf{\Sigma^{e}})_{(g,1,g^{-1})}$ has
dimension zero. So from the definition of orbifold cup product in
[5] it is easy to check that $y^{g}\cup_{orb}y^{b_{i}}=y^{g}\cdot
y^{b_{i}},~~y^{g}\cup_{orb}\gamma=y^{g}\cdot \gamma$.

For the orbifold product $y^{g_{1}}\cup_{orb}y^{g_{2}}$, where
$g_{1},g_{2}\in Box(\mathbf{\Sigma^{e}})$. From  (\ref{3-sector}),
we see that if there is no cone in $\Sigma$ containing
$\overline{g}_{1},\overline{g}_{2}$, then there is no 3-twisted
sector corresponding to the elements $g_{1},g_{2}$, so the
orbifold cup product is zero from the definition. On the other
hand  from the definition of the group ring
$A^{*}(B)[N]^{\mathbf{\Sigma^{e}}}$, $y^{g_{1}}\cdot y^{g_{2}}=0$,
so $y^{g_{1}}\cup_{orb}y^{g_{2}}=y^{g_{1}}\cdot y^{g_{2}}$. If
there is a cone $\sigma\in \Sigma$ such that
$\overline{g}_{1},\overline{g}_{2}\in \sigma$, let $g_{3}\in
Box(\mathbf{\Sigma^{e}})$ such that $\overline{g}_{3}\in
\sigma(\overline{g}_{1},\overline{g}_{2})$ and $g_{1}g_{2}g_{3}=1$
in the local group. Using the same method in [4], we get:
$y^{g_{1}}\cup_{orb}y^{g_{2}}=y^{g_{1}}\cdot y^{g_{2}}$. The
theorem is proved.  $\square$
%%% ----------------------------------------------------------------------
\section{The $\mu$-Gerbe.}

In this section we talk about the degenerate case of the extended
toric Deligne-Mumford stacks. In this case $N$ is a finite abelian
group, the simplicial fan $\Sigma$ is $0$. The toric stack bundle
is a $\mu$-gerbe $\mathcal{X}$ over $B$ for a finite abelian group
$\mu$.

Let
$N=\mathbb{Z}_{p_{1}^{n_{1}}}\oplus\cdots\oplus\mathbb{Z}_{p_{s}^{n_{s}}}$
be a finite abelian group, where $p_{1},\cdots,p_{s}$ are prime
numbers and $n_{1},\cdots,n_{s}>1$. Let $\beta^{e}:
\mathbb{Z}\longrightarrow N$ be given by the vector
$(1,1,\cdots,1)$. Because $N_{\mathbb{Q}}=0$, so $\Sigma=0$, then
$\mathbf{\Sigma^{e}}=(N,\Sigma,\beta^{e})$ is an extended stacky
fan from Section 2.1. Let
$n=lcm(p_{1}^{n_{1}},\cdots,p_{s}^{n_{s}})$, then
$n=p_{i_{1}}^{n_{i_{1}}}\cdots p_{i_{t}}^{n_{i_{t}}}$, where
$p_{i_{1}},\cdots, p_{i_{t}}$ are the distinct prime number which
have the highest powers $n_{i_{1}},\cdots, n_{i_{t}}$. Note that
the vector $(1,1,\cdots,1)$ generates an order $n$ cyclic subgroup
of $N$. We calculate the Gale dual $(\beta^{e})^{\vee}:
\mathbb{Z}\stackrel{}{\longrightarrow} \mathbb{Z}\oplus
\bigoplus_{i\notin
\{i_{1},\cdots,i_{t}\}}\mathbb{Z}_{p_{i}}^{n_{i}}$, where
$DG(\beta^{e})=\mathbb{Z}\oplus \bigoplus_{i\notin
\{i_{1},\cdots,i_{t}\}}\mathbb{Z}_{p_{i}}^{n_{i}}$, so we have the
following exact sequence:
$$0\longrightarrow \mathbb{Z}\longrightarrow \mathbb{Z}\stackrel{\beta^{e}}{\longrightarrow}N\longrightarrow
\bigoplus_{i\notin
\{i_{1},\cdots,i_{t}\}}\mathbb{Z}_{p_{i}}^{n_{i}}\longrightarrow
0$$
$$0\longrightarrow 0\longrightarrow \mathbb{Z}\stackrel{(\beta^{e})^{\vee}}{\longrightarrow}\mathbb{Z}\oplus
\bigoplus_{i\notin
\{i_{1},\cdots,i_{t}\}}\mathbb{Z}_{p_{i}}^{n_{i}}\longrightarrow
\mathbb{Z}_{n}\oplus\bigoplus_{i\notin
\{i_{1},\cdots,i_{t}\}}\mathbb{Z}_{p_{i}}^{n_{i}}\longrightarrow
0$$ So we obtain \begin{equation}\Label{gerbe}1\longrightarrow
\mu\longrightarrow \mathbb{C}^{\times}\times \prod_{i\notin
\{i_{1},\cdots,i_{t}\}}\mu_{p_{i}}^{n_{i}}\stackrel{\alpha^{e}}{\longrightarrow}\mathbb{C}^{\times}\longrightarrow
1 \end{equation} where the map $\alpha^{e}$ in (\ref{gerbe}) is
given by the matrix  $\left[
\begin{array}{c}
  n \\
  0 \\
  \vdots \\
  0\\
\end{array}
\right]$ and $\mu=\mu_{n}\times\prod_{i\notin
\{i_{1},\cdots,i_{t}\}}\mu_{p_{i}}^{n_{i}}\cong N$. The extended
toric Deligne-Mumford stack is
$\mathcal{X}(\mathbf{\Sigma^{e}})=[\mathbb{C}^{\times}/\mathbb{C}^{\times}\times
\prod_{i\notin
\{i_{1},\cdots,i_{t}\}}\mu_{p_{i}}^{n_{i}}]=\mathcal{B}\mu$, the
classifying stack of the group $\mu$. Let $L$ be a line bundle
over a smooth variety $B$, let $L^{\times}$ be  the principal
$\mathbb{C}^{\times}$-bundle induced from $L$ removing the zero
section.  From our twist we have
$^{L^{\times}}\mathcal{X}(\mathbf{\Sigma^{e}})=L^{\times}\times_{\mathbb{C}^{\times}}[\mathbb{C}^{\times}/\mathbb{C}^{\times}\times
\prod_{i\notin \{i_{1},\cdots,i_{t}\}}\mu_{p_{i}}^{n_{i}}]
=[L^{\times}/\mathbb{C}^{\times}\times \prod_{i\notin
\{i_{1},\cdots,i_{t}\}}\mu_{p_{i}}^{n_{i}}]$, which is exactly a
$\mu$-gerbe $\mathcal{X}$ over $B$. The structure of this gerbe is
a $\mu_{n}$-gerbe coming from the line bundle $L$ plus a trivial
$\prod_{i\notin \{i_{1},\cdots,i_{t}\}}\mu_{p_{i}}^{n_{i}}$-gerbe
over $B$. For this toric stack bundle, we know that 
$Box(\mathbf{\Sigma^{e}})=N$, so we have the following Proposition
for the inertia stack.

\begin{prop}
The inertia stack of this toric stack bundle $\mathcal{X}$ is
$p_{1}^{n_{1}}\cdot\cdots\cdot p_{s}^{n_{s}}$ copies of the
$\mu$-gerbe $\mathcal{X}$.
\end{prop}
From our main Theorem, we have:

\begin{prop}
The orbifold cohomology ring of the toric stack bundle
$\mathcal{X}$ is given by:
$$H^{*}_{orb}(\mathcal{X},\mathbb{Q})\cong H^{*}(B,\mathbb{Q})\otimes H^{*}_{orb}(\mathcal{B}\mu,\mathbb{Q})$$
where
$H^{*}_{orb}(\mathcal{B}\mu;\mathbb{Q})=\mathbb{Q}[t_{1},\cdots,t_{s}]/(t_{i}^{p_{i}^{n_{i}}}-1)$.
\end{prop}

Let $N=\mathbb{Z}_{r}$, and $\beta: \mathbb{Z}\longrightarrow
\mathbb{Z}_{r}$ be the natural projection. The toric
Deligne-Mumford stack
$\mathcal{X}(\mathbf{\Sigma})=\mathcal{B}\mu_{r}$. Let
$L\longrightarrow B$ be a line bundle, then the toric stack bundle
$\mathcal{X}=B_{(L,r)}$ is the $\mu_{r}$-gerbe over $B$ determined
by the line bundle $L$. We have:

\begin{cor}
The orbifold cohomology ring of $B_{(L,r)}$ is isomorphic to
$H^{*}(B)[t]/(t^{r}-1)$.
\end{cor}

If the variety $B$ is not a toric variety, then the toric stack
bundle over $B$ is not a toric Deligne-Mumford stack. But suppose
$B$ is a smooth toric variety, then a $\mu$-gerbe $\mathcal{X}$
can give a toric Deligne-Mumford stack in [4].

\begin{example}
Let $B=\mathbb{P}^{d}$ be the $d$-dimensional projective space. We
give stacky fan $\mathbf{\Sigma}=(N,\Sigma,\beta)$ as follows: let
$N=\mathbb{Z}^{d}\oplus \mathbb{Z}_{r}$, $\beta:
\mathbb{Z}^{d+1}\longrightarrow N$ be determined by the vectors:
$\{(1,0,\ldots,0,0), (0,1,\ldots,0,0),\ldots, (0,0,\ldots,1,0),
(-1,-1,\ldots,-1,1)\}$. Then $DG(\beta)=\mathbb{Z}$, and the Gale
dual $\beta^{\vee}$ is given by the matrix $[r,r,\ldots,r]$. So we
have the following exact sequences:
$$0\longrightarrow \mathbb{Z}\longrightarrow \mathbb{Z}^{d+1}\stackrel{\beta}
{\longrightarrow}\mathbb{Z}^{d}\oplus\mathbb{Z}_{r}\longrightarrow 0\longrightarrow 0$$
$$0\longrightarrow \mathbb{Z}^{d}\longrightarrow \mathbb{Z}^{d+1}\stackrel
{\beta^{\vee}}{\longrightarrow}\mathbb{Z}\longrightarrow \mathbb{Z}_{r}\longrightarrow 0$$
Then we obtain the exact sequence:
$$1\longrightarrow \mu_{r}\longrightarrow \mathbb{C}^{\times}\stackrel{\alpha}
{\longrightarrow}(\mathbb{C}^{\times})^{d+1}\longrightarrow (\mathbb{C}^{\times})^{d}\longrightarrow 1$$
The toric Deligne-Mumford stack
$\mathcal{X}(\mathbf{\Sigma}):=[\mathbb{C}^{d+1}-\{0\}/\mathbb{C}^{\times}]$
is the canonical $\mu_{r}$-gerbe over the projective space
$\mathbb{P}^{d}$ coming from the canonical line bundle, where the
$\mathbb{C}^{\times}$ action is given by $\lambda\cdot
(z_{1},\ldots,z_{d+1})=(\lambda^{r}\cdot
z_{1},\ldots,\lambda^{r}\cdot z_{d+1})$. Denote this toric
Deligne-Mumford stack by $\mathcal{G}_{r}=\mathbb{P}(r,\ldots,r)$.
If the homomorphism $\beta: \mathbb{Z}^{d+1}\longrightarrow N$ is
determined by the vectors: $\{(1,0,\ldots,0,0),
(0,1,\ldots,0,0),\ldots, (0,0,\ldots,1,0), (-1,-1,\ldots,-1,0)\}$,
then $DG(\beta)=\mathbb{Z}\oplus \mathbb{Z}_{r}$. Comparing to the
former exact sequence, we have the exact sequence:
$$1\longrightarrow \mu_{r}\longrightarrow \mathbb{C}^{\times}\times
\mu_{r}\stackrel{\alpha}{\longrightarrow}(\mathbb{C}^{\times})^{d+1}\longrightarrow
(\mathbb{C}^{\times})^{d}\longrightarrow 1$$ The corresponding
toric Deligne-Mumford stack is the trivial $\mu_{r}$-gerbe
$\mathbb{P}^{d}\times \mathcal{B}\mu_{r}$ coming from the trivial
line bundle over $\mathbb{P}^{d}$. The coarse Moduli spaces of
these two stacks are all the projective space $\mathbb{P}^{d}$.
From the Theorem of this paper or the main Theorem in [4], the
orbifold cohomology rings of these two stacks are isomorphic
although as stacks, they are different.
\end{example}

\begin{rmk}
Let $H$ represent the hyperplane class of  $\mathbb{P}^{d}$, then
$H^{*}_{orb}(\mathcal{G}_{r},\mathbb{Q})\cong
\mathbb{Q}[H]/(H^{d+1})\otimes \mathbb{Q}[t]/(t^{r}-1)$. We
conjecture that the orbifold quantum cohomology ring of
$\mathcal{G}_{r}$ defined in [6] is isomorphic to
$\mathbb{Q}[H]/(H^{d+1}-f(H,q))\otimes
\mathbb{Q}[t]/(t^{r}-1-g(t,q))$, where $f,g$ are two relations and
$q$ is the quantum parameter. The orbifold quantum cohomology of
trivial gerbe case has been computed in [12], where $f(H,q)=q$ and
$g(t,q)=0$.
\end{rmk}

\begin{rmk}
We conjecture that the small orbifold quantum cohomology ring of
the nontrivial $\mu_{r}$-gerbe and trivial $\mu_{r}$-gerbe over
the projective space $\mathbb{P}^{d}$ should be different. This
means that the orbifold quantum cohomology can classify these two
different stacks.
\end{rmk}
%%% ----------------------------------------------------------------------
\section{Application.}
In this section we generalize one result of Borisov, Chen and
Smith [4] to the toric stack bundle case.

Let $X(\Sigma)$ be a simplicial toric variety, and let
$\mathcal{X}(\mathbf{\Sigma})$ be the associated toric DM stack,
where $\mathbf{\Sigma}=(N,\Sigma,\beta)$.  Let $\Sigma^{'}$ be a
subdivision of $\Sigma$ such that $X(\Sigma^{'})$ is a crepant
resolution of $X(\Sigma)$. Suppose there are $m$ rays in
$\Sigma^{'}$, let $i: (\mathbb{C}^{\times})^{n}\longrightarrow
(\mathbb{C}^{\times})^{m}$ be the inclusion. From the following
commutative diagram:
\[
\begin{CD}
0 @ >>>\mathbb{Z}^{n-d}@ >>> \mathbb{Z}^{n}@ >{\beta}>> N @
>>> 0\\
&& @VV{}V@VV{i}V@VV{id}V \\
0@ >>>\mathbb{Z}^{m-d} @ >{}>>\mathbb{Z}^{m}@ >\beta^{'}>> N @>>>
0
\end{CD}
\]
Taking Gale dual we get:
\[
\begin{CD}
0 @ >>>N^{*}@ >>> (\mathbb{Z}^{m})^{*}@ >{(\beta^{'})^{\vee}}>>
DG(\beta^{'}) @
>>> 0\\
&& @VV{id}V@VV{}V@VV{}V \\
0@ >>>N^{*} @ >{}>>(\mathbb{Z}^{n})^{*}@ >{\beta^{\vee}}>>
DG(\beta) @>>> 0
\end{CD}
\]
So applying the $Hom$ functor we have the following diagram:
$$\xymatrix{
(\mathbb{C}^{\times})^{n}\dto^{}\rto^{i} &(\mathbb{C}^{\times})^{m}\dto^{}\\
T\rto^{id} &T }$$ Let $P\longrightarrow B$ be a principal
$(\mathbb{C}^{\times})^{n}$-bundle, we still use $P$ to represent
the principal $(\mathbb{C}^{\times})^{m}$-bundle induced by $i$,
then they induce the same principal $T$ bundle $E$ over $B$. So
$^{E}X(\Sigma^{'})\longrightarrow ^{E}X(\Sigma)$ is the crepant
resolution. And $^{E}X(\Sigma)$ is the coarse moduli space of the
toric stack bundle $^{P}\mathcal{X}(\mathbf{\Sigma})$ from
Proposition 3.3. We have the following result:

\begin{prop}
If the Chow ring of the smooth variety $B$ is a Cohen-Macaulay
ring. Then there is a flat family $\mathcal{S}\longrightarrow
\mathbb{P}^{1}$ of schemes such that $\mathcal{S}_{0}\cong
Spec(A^{*}_{orb}(^{P}\mathcal{X}(\mathbf{\Sigma})))$ and
$\mathcal{S}_{\infty}\cong Spec(A^{*}(^{E}X(\Sigma^{'})))$.
\end{prop}
\begin{pf}
We also construct a family of algebras over $\mathbb{P}^{1}$ such
that the fiber over $0$ and $\infty$ are
$A^{*}_{orb}(^{P}\mathcal{X}(\mathbf{\Sigma}))$ and
$A^{*}_{orb}(^{P}\mathcal{X}(\mathbf{\Sigma}))$ respectively.
Since $X(\Sigma^{'})$ is a smooth variety,
$\{b_{1},\cdots,b_{n},b_{n+1},\cdots b_{m}\}$ generate the whole
lattice $N$, then $A^{*}(B)[N]^{\mathbf{\Sigma}}$ is the quotient
ring of the ring $S=A^{*}(B)[y^{b_{1}},\cdots,y^{b_{m}}]$ by the
binomial ideal determined by (\ref{product}). Let $I_{2}$ denote
this ideal. Let $I_{1}$ denote the ideal generated by
$c_{1}(\xi_{\theta_{j}})+\sum_{i=1}^{m}\theta_{j}(b_{i})y^{b_{i}}$
for $1\leq j\leq d$, where $\theta_{1},\cdots,\theta_{d}$ is a
basis of $N^{*}$. Since $\Sigma^{'}$ is a regular subdivision of
$\Sigma$, then there is a $\Sigma^{'}$-linear support function $h:
N\longrightarrow \mathbb{Z}$ such that $h(b_{i})=0$ for $1\leq
i\leq n$, $h(b_{i})>0$ for $n+1\leq i\leq m$. For any lattice
points $c_{1},c_{2}$ lying in the same cone of $\Sigma$,
$h(c_{1}+c_{2})\geq h(c_{1})+h(c_{2})$, and the inequality is
strict unless $c_{1},c_{2}$ lies  in the same cone of
$\Sigma^{'}$.

Let $\breve{I}_{1}$ be the ideal in $S[t_{1}]$ generated by
$c_{1}(\xi_{\theta_{j}})t_{1}^{h(b_{i})}+\sum_{i=1}^{m}\theta_{j}(b_{i})y^{b_{i}}t_{1}^{h(b_{i})}$
for $1\leq j\leq d$. So the choice of $h$ implies that
$$\frac{S[t_{1}]}{\breve{I}_{1}+I_{2}+<t_{1}>}\cong \frac{S}
{<c_{1}(\xi_{\theta_{j}})+\sum_{i=1}^{n}\theta_{j}(b_{i})y^{b_{i}}:1\leq
j\leq d>+I_{2}}\cong A^{*}(^{P}\mathcal{X}(\mathbf{\Sigma}))$$ The
sequence
$c_{1}(\xi_{\theta_{j}})+\sum_{i=1}^{n}\theta_{j}(b_{i})y^{b_{i}}$
for $1\leq j\leq d$ is also a homogeneous system of parameters on
$S/I_{2}$. The Chow ring  $A^{*}(B)$ is a Cohen-Macaulay ring, so
$S/I_{2}$ is also Cohen-Macaulay. So the sequence is a regular
sequence. Therefore, the Hilbert function of the family
$S[t_{1}/(\breve{I}_{1}+I_{2})]$ is constant outside a finite set
in $\mathbb{Q}^{*}$.

On the other hand, for the family over $\mathbb{P}^{1}-\{0\}$, we
use the similar method in [4] to get the family
$S[t_{2}]/(I_{1}+\breve{I}_{2})$, where $\breve{I}_{2}$ is the
binomial ideal in $S[t_{2}]$ defined in [4]. There exists an
automorphism $\varphi$ betweem these two families so that we
construct such a family over $\mathbb{P}^{1}$. All the left proof
is the same as the one in [4]. We omit the details.
\end{pf}

\begin{rmk}
Ruan [18] conjectured that the cohomology ring of crepant
resolution is isomorphic to the orbifold Chow ring of the orbifold
if we add some quantum corrections on the ordinary cohomology ring
of the crepant resolution coming from the exceptional divisors.
Let $\mathbb{P}(1,1,2)$ be the weighted projective plane with one
orbifold point with local group $\mathbb{Z}_{2}$, the Hirzburch
surface is the crepant resolution of $\mathbb{P}(1,1,2)$.  We can
compute the quantum correction of the cohomology ring of the
Hirzburch surface and check Ruan's conjecture. This case has been
done recently in [16].
\end{rmk}

% ------------------------------------------------------------------------

\subsection*{}

% ------------------------------------------------------------------------
\author{The University of British Columbia, Vancouver, BC,
Canada}\\
jiangyf@math.ubc.ca
\end{document}